\documentclass[11pt,twoside,final]{scrartcl}
\usepackage{url,a4wide}
\usepackage{amsmath, amsfonts, amssymb, amsthm, mathtools, nicefrac}
\usepackage{todonotes}
\usepackage{enumitem}
\usepackage{graphicx}
\graphicspath{{pics/}}
\usepackage[skip=1pt, font=normalsize]{subcaption}
\usepackage{hyperref} % provides \url command for bibtex and links to jump within documents
\hypersetup{plainpages=false, colorlinks, linkcolor=black, citecolor=black, urlcolor=blue,
pdftitle={On regularized Shannon sampling formulas with localized sampling},
pdfauthor={Melanie Kircheis, Daniel Potts, Manfred Tasche},
pdfstartview={FitBH}}
\allowdisplaybreaks

\newcommand{\e}{\mathrm e}

\newcommand{\sinc}{\mathrm{sinc}}

\newcommand{\R}{\mathbb R}
\newcommand{\C}{\mathbb C}
\newcommand{\Z}{\mathbb Z}
\newcommand{\N}{\mathbb N}

\renewcommand{\O}[1]{\mathcal O(#1)}

\newtheorem{theorem}{Theorem}[section]
\newtheorem{corollary}[theorem]{Corollary}
\newtheorem{lemma}[theorem]{Lemma}
\newtheorem{alg}[theorem]{Algorithm}
\newtheorem{definition}[theorem]{Definition}
\newtheorem{example}[theorem]{Example}
\newtheorem{remark}[theorem]{Remark}

\newcommand{\bend}{\hspace*{0ex} \hfill \hbox{\vrule height
	1.5ex\vbox{\hrule width 1.4ex \vskip 1.4ex\hrule  width 1.4ex}\vrule
	height 1.5ex}}
\renewcommand{\qedsymbol}{\rule{1.5ex}{1.5ex}}

\newenvironment{Lemma}{\goodbreak\begin{lemma}}{\end{lemma}}
\newenvironment{Theorem}{\goodbreak\begin{theorem}}{\end{theorem}}
\newenvironment{Remark}{\goodbreak\begin{remark}\upshape}{\bend\end{remark}}

\newenvironment{Example}{\goodbreak\begin{example}\upshape}{\bend\end{example}}
\newenvironment{algorithm}[1]{\goodbreak~\begin{alg}[#1]~\vspace{-9pt}~\\
		\rule{\linewidth}{0.5pt}~\\}{\vspace{-9pt}~\\
		\rule{\linewidth}{0.5pt}~\end{alg}}

\numberwithin{equation}{section}
\numberwithin{table}{section}
\numberwithin{figure}{section}

\renewcommand{\mathbf}[1]{\ensuremath{\boldsymbol{#1}}}

\title{On regularized Shannon sampling formulas with localized sampling}
\author{
	Melanie Kircheis\footnotemark[1] \and
	Daniel Potts\footnotemark[4] \and
	Manfred Tasche\footnotemark[3]
}

\date{}
\begin{document}
\maketitle

\begin{abstract}
	In this paper we present new regularized Shannon sampling formulas which use localized sampling with special window functions,
	namely Gaussian, $\mathrm B$--spline, and $\sinh$-type window functions. In contrast to the classical Shannon sampling series, the regularized Shannon sampling formulas possess an
	exponential decay and are numerically robust in the presence of noise. Several numerical experiments illustrate the theoretical results.
	\medskip
	
	\smallskip
	
	AMS \emph{Subject Classifications}:
	94A20, 65T50.
\end{abstract}

\footnotetext[1]{Corresponding author: melanie.kircheis@math.tu-chemnitz.de, Chemnitz University of Technology, Faculty of Mathematics, D--09107 Chemnitz, Germany}
\footnotetext[4]{potts@mathematik.tu-chemnitz.de, Chemnitz University of
	Technology, Faculty of Mathematics, D--09107 Chemnitz, Germany}
\footnotetext[3]{manfred.tasche@uni-rostock.de, University of Rostock, Institute of Mathematics, D--18051 Rostock, Germany}

\section{Introduction}

The classical Whittaker--Kotelnikov--Shannon sampling theorem plays a fundamental role in signal processing, since it describes the close relation between a bandlimited
function and its equidistant samples.
More precisely, this sampling theorem states that any function $f \in L^2(\mathbb R)$ with bandwidth $\le \frac{N}{2}$, i.\,e.,
the support of the Fourier transform
$$
{\hat f}(v) \coloneqq \int_{\mathbb R} f(t)\,{\mathrm e}^{-2 \pi {\mathrm i}t v}\,{\mathrm d}t\,, \quad v \in \mathbb R\,,
$$
is contained in $\big[- \frac{N}{2},\,\frac{N}{2}\big]$, can be recovered from its samples $f\big(\frac{\ell}{L}\big)$, $\ell \in \mathbb Z$, with $L \ge N$
and it holds
\begin{equation}
\label{eq:WKSseries}
f(t) = \sum_{\ell \in \mathbb Z} f\big(\tfrac{\ell}{L}\big)\,\mathrm{sinc}\big(L\pi\,\big(t - \tfrac{\ell}{L}\big)\big)\,, \quad t \in \mathbb R\,,
\end{equation}
with the $\sinc$ function
$$
\mathrm{sinc}\,x \coloneqq \left\{ \begin{array}{ll}  \frac{\sin x}{x} & \quad x \in \mathbb R \setminus \{0\}\,, \\ [1ex]
1 & \quad x = 0\,. \end{array} \right.
$$
Unfortunately, the practical use of this sampling theorem is limited, since it requires infinitely many samples, which is impossible in practice.
Further, the $\sinc$ function decays very slowly such that the \emph{Shannon sampling series}
$$
\sum_{\ell \in \mathbb Z} f\big(\tfrac{\ell}{L}\big)\,\mathrm{sinc}\big(L\pi\,\big(t - \tfrac{\ell}{L}\big)\big)\,, \quad t \in \mathbb R\,,
$$
has rather poor convergence.
Moreover, in the presence of noise or quantization in the samples $f\big(\frac{\ell}{L}\big)$, $\ell \in \mathbb Z$, the convergence of Shannon
sampling series may even break down completely (see \cite{DDeV03}).

To overcome these drawbacks, one can use the following three techniques (see \cite{Q03, LZ16, MXZ09, StTa06}):\\
1. The function ${\mathrm{sinc}}(L \pi\, \cdot)$ is \emph{regularized by a truncated window function}
$$
\varphi_m(x) \coloneqq \varphi(x)\,{\mathbf 1}_{[-m/L,\,m/L]}(x)\,, \quad x \in \mathbb R\,,
$$
where the window function $\varphi: \mathbb R \to [0,\,1]$ belongs to the set $\Phi_{m,L}$ (as defined in Section 3) and ${\mathbf 1}_{[-m/L,\,m/L]}$ denotes
the indicator function of the
interval $\big[-\frac{m}{L},\,\frac{m}{L}\big]$ with some \mbox{$m \in \mathbb N \setminus \{1\}$}. Then we recover a function $f \in L^2(\mathbb R)$ with
bandwidth $\le \frac{N}{2}$ by the \emph{regularized Shannon sampling formula}
\begin{align*}
	(R_{\varphi,m} f)(t) \coloneqq \sum_{\ell \in \mathbb Z} f\left(\tfrac{\ell}{L}\right)\, {\mathrm{sinc}}\big(L \pi\left(t - \tfrac{\ell}{L}\right)\big)\,\varphi_m\left(t - \tfrac{\ell}{L}\right) \,,
\end{align*}
where $L \ge N$. Obviously,
this is an \emph{interpolating approximation} of $f$, since it holds
$$
(R_{\varphi,m} f)\left(\tfrac{k}{L}\right) = f\left(\tfrac{k}{L}\right)\,, \quad k \in \mathbb Z\,.
$$
2. The use of the truncated window function $\varphi_m$ with compact support $\big[- \frac{m}{L},\,\frac{m}{L}\big]$ leads to \emph{localized sampling} of $f$, i.\,e., the
computation of $(R_{\varphi,m} f)(t)$ for $t \in \mathbb R \setminus \frac{1}{L}\,\mathbb Z$ requires only $2m$ samples $f\big(\frac{k}{L}\big)$, where $k\in \mathbb Z$ fulfills the
condition $|k - L t| \le m$.
If $f$ has bandwidth $\le N/2$ and $L \ge N$, then the reconstruction of $f$ on the interval $[0,\,1]$ requires only $2m + L$
samples $f\big(\frac{\ell}{L}\big)$
with $\ell = -m,\,1-m,\,\ldots,\,m+L$.\\
3. In many applications, one usually employs \emph{oversampling}, i.\,e., a function $f \in L^2(\mathbb R)$ of bandwidth $\le N/2$ is sampled on a finer grid $\frac{1}{L}\, \mathbb Z$ with $L > N$.

This concept of regularized Shannon sampling formulas with localized sampling and oversampling has already been studied by various authors, e.\,g. in \cite{Q03, Q06} and references therein for the Gaussian window function.
An improvement of the theoretical error bounds for the Gaussian window function was made by \cite{LZ16}, whereas oversampling was neglected in this work. Rather, the special case $L=N=1$ was studied.
The case of erroneous sampling for the Gaussian window function was examined in \cite{Q05}.
Generalizations of the Gaussian regularized Shannon sampling formula to holomorphic functions $f$ were introduced by \cite{SchSt07} using contour integration and by \cite{TaSuMu07} for the approximation of derivatives of $f$.
A survey of different approaches for window functions can be found in \cite{Q04}.
Furthermore, in \cite{StTa06} the problem was approached in Fourier space.
Oversampling then is equivalent to continuing the Fourier transform of the $\sinc$ function on the larger interval $[- L/2, L/2]$.
Here the aim is to find a regularization function whose Fourier transform is smooth.
However, the resulting function does not have an explicit representation and therefore cannot directly be used in spatial domain.
Nevertheless, the complexity and efficiency of the received methods was not the main focus of the aforementioned approaches.
On the contrary, we now propose a new set $\Phi_{m,L}$ of window functions $\varphi$ such that smaller truncation parameters $m$ are sufficient for achieving high accuracy, therefore yielding short sums being evaluable very fast.

In this paper we present new regularized Shannon sampling formulas with localized sampling and derive new estimates of the uniform approximation error
$$
\max_{t \in \mathbb R}\big| f(t) - (R_{\varphi,m} f)(t) \big|\,,
$$
where we apply several window functions $\varphi$.
It is shown in the subsequent sections that the uniform approximation error decays exponentially with respect to $m$, if \mbox{$\varphi \in \Phi_{m,L}$} is the Gaussian, $\mathrm B$--spline, or $\sinh$-type window function.
Otherwise, if \mbox{$\varphi \in \Phi_{m,L}$} is chosen as the rectangular window function, then the uniform approximation error of the regularized Shannon sampling formula has a poor decay of order $m^{-1/2}$.
Moreover, we show that the regularized Shannon sampling formulas are numerically robust for noisy samples, i.\,e., if \mbox{$\varphi \in \Phi_{m,L}$} is the Gaussian, $\mathrm B$--spline, or $\sinh$-type window function,
then the uniform perturbation error \eqref{eq:perturb_error} only grows as $m^{1/2}$.

In our approach we need the Fourier transform of the product of ${\mathrm{sinc}}(L \pi\, \cdot)$ and the window function $\varphi$.
Since the ${\mathrm{sinc}}$ function belongs to $L^2(\mathbb R)$, but not to $L^1(\mathbb R)$, we present the convolution property of the Fourier transform for $L^2(\mathbb R)$ functions in the preliminary Section~2 (see Lemma 2.2).
In Section~\ref{sec:WKS_localized} we consider regularized Shannon sampling formulas for an arbitrary window function $\varphi \in \Phi_{m,L}$.
Here the main results are Theorem~\ref{Theorem:errorRmf} with a unified approach to error estimates for regularized Shannon sampling formulas and Theorem~\ref{Theorem:robustness} on the numerical robustness of regularized Shannon sampling formulas.
Afterwards, we concretize the results for several choices of window functions.
In Section~\ref{sec:Gaussian} we consider the Gaussian window function (as in \cite{Q03, LZ16}).
Then Theorem~\ref{Theorem:GaussianWKS} shows that the uniform approximation error decays exponentially with respect to $m$.
In Section~\ref{sec:Bspline} we use the $\mathrm B$--spline window function and prove the exponential decay of the uniform approximation error with respect to $m$ in Theorem~\ref{Theorem:bsplineWKS}.
In Section~\ref{sec:sinh} we discuss the $\sinh$-type window function, where it is proven in Theorem~\ref{Thm:errorpsisinh} that the uniform approximation error decays exponentially with respect to $m$. Several numerical experiments illustrate the theoretical results.
Finally, in the concluding Section~\ref{sec:Conclusion}, we compare the proposed window functions and show the superiority of the new $\sinh$-type window function.

\section{Convolution property of the Fourier transform}

Let $C_0(\mathbb R)$ denote the Banach space of continuous functions $f \colon \R\to\C$ vanishing as $|x|\to\infty$ with norm
$$
\|{f}\|_{C_0(\mathbb R)} \coloneqq \max_{t\in \mathbb R} |{f}(t)|\,.
$$
As known (see \cite[p.~66]{PPST18}), the Fourier transform defined by
\begin{equation}
\label{eq:Fouriertrafo}
	{\hat f}(v) \coloneqq \int_{\mathbb R} f(t)\,{\mathrm e}^{-2 \pi {\mathrm i}v t}\, {\mathrm d}t\,, \quad v \in \mathbb R\,,
\end{equation}
is a continuous mapping from $L^1(\mathbb R)$ into $C_0(\mathbb R)$ with
$$\|\hat f\|_{C_0(\mathbb R)} \le \| f\|_{L^1(\mathbb R)} \coloneqq \int_{\mathbb R} |f(t)|\, {\mathrm d}t\,.
$$
Here we are interested in the Hilbert space $L^2(\mathbb R)$ with inner product and norm
\begin{align*}
	\langle f , g \rangle_{L^2(\mathbb R)} \coloneqq \int_{\mathbb R} f(t) \, \overline{g(t)} \, {\mathrm d}t \,, \quad
	\|f\|_{L^2(\mathbb R)} \coloneqq \left( \langle f , f \rangle_{L^2(\mathbb R)}\right)^{1/2} \,.
\end{align*}
By the theorem of Plancherel, the Fourier transform is also an invertible, continuous mapping from $L^2(\mathbb R)$ onto itself with $\|f\|_{L^2(\mathbb R)} = \|\hat f\|_{L^2(\mathbb R)}$.

For $f,\,g \in L^1(\mathbb R)$ the convolution property of the Fourier transform reads as
\begin{equation}
	\label{eq:convL1}
	(f * g){\hat{}} = {\hat f}\,{\hat g} \in C_0(\mathbb R)\,,
\end{equation}
where the convolution is defined by
$$
(f * g)(x) \coloneqq \int_{\mathbb R} f(x-t)\,g(t)\, {\mathrm d}t\,, \quad x \in \mathbb R\,.
$$
However, for any $f$, $g \in L^2(\mathbb R)$ the convolution property of the Fourier transform is not true in the form \eqref{eq:convL1}, since by Young's inequality $f * g \in C_0(\mathbb R)$, by H\"{o}lder's inequality
${\hat f}\,{\hat g} \in L^1(\mathbb R)$ and since the Fourier transform does not map $C_0(\mathbb R)$ into $L^1(\mathbb R)$.
Instead, the convolution property
of the Fourier transform in $L^2(\mathbb R)$ has the following form.

\begin{Lemma}
	\label{Lemma:convL2}
	For all $f,\,g \in L^2(\mathbb R)$ it holds
	\begin{equation}
		\label{eq:convL2}
		f * g = ({\hat f}\, {\hat g}){\check{}} \in C_0(\mathbb R)\,,
	\end{equation}
	where $\check h$ denotes the inverse Fourier transform of $h \in L^1(\mathbb R)$ defined as
	\begin{equation}
		\label{eq:inverse_Fouriertrafo}
		{\check h}(t) \coloneqq \int_{\mathbb R} h(v)\,{\mathrm e}^{2 \pi {\mathrm i}v t}\,{\mathrm d}v\,, \quad t \in \mathbb R\,.
	\end{equation}
\end{Lemma}

\emph{Proof}. For arbitrary $f,\,g \in L^2(\mathbb R)$ it holds ${\hat f},\,{\hat g} \in L^2(\mathbb R)$, see \cite[p.~80]{PPST18}.
Since the Schwartz space ${\mathcal S}(\mathbb R)$, cf.~\cite[p.~167]{PPST18}, is dense in $L^2(\mathbb R)$, see.~\cite[p.~170]{PPST18}, there exist sequences $(f_n)_{n=1}^{\infty}$ and $(g_n)_{n=1}^{\infty}$ in ${\mathcal S}(\mathbb R)$ such that
\begin{equation}
	\label{eq:limfngn}
	\lim_{n\to \infty} \| f_n - f \|_{L^2(\mathbb R)} = \lim_{n\to \infty} \| g_n - g \|_{L^2(\mathbb R)} = 0\,.
\end{equation}
Since the Fourier transform is a continuous mapping on $L^2(\mathbb R)$, it follows that
\begin{equation}
	\label{eq:limhatfnhatgn}
	\lim_{n\to \infty} \| {\hat f}_n - {\hat f} \|_{L^2(\mathbb R)} = \lim_{n\to \infty} \| {\hat g}_n - {\hat g} \|_{L^2(\mathbb R)} = 0\,.
\end{equation}
If we now write
$$
(f * g) - (f_n * g_n) = (f - f_n) * g + f_n * (g - g_n)\,,
$$
we see by the triangle inequality and Young's inequality that
$$
\|(f * g) - (f_n * g_n)\|_{C_0(\mathbb R)} \le \|f - f_n\|_{L^2(\mathbb R)}\,\| g \|_{L^2(\mathbb R)}+ \| f_n \|_{L^2(\mathbb R)}\, \| g - g_n\|_{L^2(\mathbb R)}
$$
and hence by \eqref{eq:limfngn} that
\begin{equation}
	\label{eq:limfnastgn}
	\lim_{n \to \infty} \| (f * g) - (f_n * g_n)\|_{C_0(\mathbb R)} = 0\,.
\end{equation}
On the other hand, if we write
$$
{\hat f}\,{\hat g} - {\hat f}_n\, {\hat g}_n = ({\hat f} - {\hat f}_n)\,{\hat g} + {\hat f}_n\,({\hat g} - {\hat g}_n)\,,
$$
we see by the triangle inequality and H\"{o}lder's inequality that
$$
\|{\hat f}\,{\hat g} - {\hat f}_n\, {\hat g}_n \|_{L^1(\mathbb R)} \le \|{\hat f} - {\hat f}_n\|_{L^2(\mathbb R)}\,\|{\hat g}\|_{L^2(\mathbb R)} + \|{\hat f}_n\|_{L^2(\mathbb R)}\,\|{\hat g} - {\hat g}_n\|_{L^2(\mathbb R)}
$$
and hence by \eqref{eq:limhatfnhatgn} that
\begin{equation}
	\label{eq:limhatfhatg}
	\lim_{n\to \infty} \|{\hat f}\,{\hat g} - {\hat f}_n\, {\hat g}_n \|_{L^1(\mathbb R)} = 0\,.
\end{equation}
By the convolution property of the Fourier transform in ${\mathcal S}(\mathbb R)$, we have for $f_n$, $g_n \in {\mathcal S}(\mathbb R)$ that
$$
(f_n * g_n){\hat{}} = {\hat f}_n\, {\hat g}_n\,.
$$
Note that $f_n * g_n \in {\mathcal S}(\mathbb R)$ and ${\hat f}_n\, {\hat g}_n \in {\mathcal S}(\mathbb R)$ (see \cite[p.~175]{PPST18}).
Since the Fourier transform on ${\mathcal S}(\mathbb R)$ is invertible (see \cite[p.~175]{PPST18}), it holds
\begin{equation}
	\label{eq:fnastgn}
	f_n * g_n = ({\hat f}_n \,{\hat g}_n){\check{}}\,.
\end{equation}
Moreover, since the inverse Fourier transform is a continuous mapping from $L^1(\mathbb R)$ into $C_0(\mathbb R)$, it holds by \cite[pp.~66--67]{PPST18} that
$$
\| ({\hat f}\,{\hat g}){\check{}} - ({\hat f}_n \,{\hat g}_n){\check{}}\,\|_{C_0(\mathbb R)} \le  \|{\hat f}\,{\hat g} - {\hat f}_n\, {\hat g}_n \|_{L^1(\mathbb R)}\,.
$$
From \eqref{eq:limhatfhatg} it follows that
\begin{equation}
	\label{eq:lim(hatfhatg)check}
	\lim_{n\to \infty} \| ({\hat f}\,{\hat g}){\check{}} - ({\hat f}_n \,{\hat g}_n){\check{}}\,\|_{C_0(\mathbb R)}  = 0\,.
\end{equation}
Thus, by \eqref{eq:fnastgn} we conclude that
\begin{eqnarray*}
	\|f \ast g - ({\hat f} \,{\hat g}){\check{}}\, \|_{C_0(\mathbb R)} &\le& \|f\ast g - f_n \ast g_n\|_{C_0(\mathbb R)}
	+ \|f_n \ast g_n - ({\hat f} \,{\hat g}){\check{}}\,\|_{C_0(\mathbb R)}\\
	&=& \|f\ast g - f_n \ast g_n\|_{C_0(\mathbb R)} + \|({\hat f}_n \,{\hat g}_n){\check{}}   - ({\hat f} \,{\hat g}){\check{}}\,\|_{C_0(\mathbb R)}\,.
\end{eqnarray*}
For $n\to \infty$ the right hand side of above estimate converges to zero by \eqref{eq:limfnastgn} and \eqref{eq:lim(hatfhatg)check}. This implies \eqref{eq:convL2}. \qedsymbol
\medskip

The following equivalent formulation of the convolution property in $L^2(\mathbb R)$ can be obtained, if we replace $f \in L^2(\mathbb R)$ by $\hat f \in L^2(\mathbb R)$
and $g \in L^2(\mathbb R)$ by $\hat g \in L^2(\mathbb R)$ in \eqref{eq:convL2}.

\begin{Lemma}
\label{Lemma:FTconvolution}
	For all $f$, $g \in L^2(\mathbb R)$ it holds
	\begin{equation*}
		{\hat f} * {\hat g} = (f \,g)\hat{} \in C_0(\mathbb R)\,.
	\end{equation*}
\end{Lemma}

\emph{Proof.} For any $f$, $g \in L^2(\mathbb R)$ it holds
$\hat{\hat f} = f(-\,\cdot)$ as well as $\hat{\hat g} =  g(-\, \cdot)\,.
$
Note that by H\"{o}lder's inequality we have $f\,g \in L^1(\mathbb R)$. Then by above Lemma \ref{Lemma:convL2}
it follows that
\begin{eqnarray*}
	({\hat f} * {\hat g})(t) &=& \big(\hat{\hat f}\, \hat{\hat g}\big)\check{}(t) = \big(f(-\,\cdot)\, g(-\,\cdot)\big)\check{}(t)\\
	&=& \int_{\mathbb R} f(-v)\,g(-v)\,{\mathrm e}^{2 \pi {\mathrm i}v t}\, {\mathrm d}v= \int_{\mathbb R} f(v)\,g(v)\, {\mathrm e}^{-2 \pi {\mathrm i}v t}\,{\mathrm d}v
	= (f\, g)\hat{}\,(t)\,.
\end{eqnarray*}
This completes the proof. \qedsymbol
\medskip

Note that Lemma \ref{Lemma:FTconvolution} improves a corresponding result in \cite[p.~209]{Const16}. There it was remarked that for $f$, $g \in L^2(\mathbb R)$
it holds $(f \,g)\hat{} = {\hat f} * {\hat g} \in L^{\infty}(\mathbb R)$, but by Lemma \ref{Lemma:FTconvolution} the function ${\hat f} * {\hat g}$ indeed belongs to $C_0(\mathbb R) \subset L^{\infty}(\mathbb R)$.

\section{Regularized Shannon sampling formulas with localized sampling\label{sec:WKS_localized}}

Let
$$
{\mathcal B}_{\delta}(\mathbb R) \coloneqq \{f \in L^2(\mathbb R):\, \mathrm{supp}\, \hat f
\subseteq [- \delta, \, \delta]\}
$$
be the \emph{Paley--Wiener space}.
The functions of ${\mathcal B}_{\delta}(\mathbb R)$ are called \emph{bandlimited to} $[- \delta,\,\delta]$, where $\delta > 0$ is the so-called \emph{bandwidth}.
By definition the Paley--Wiener space ${\mathcal B}_{\delta}(\mathbb R)$ consists of equivalence classes of almost equal functions.
Each of these equivalence classes contains a smooth function, since by inverse
Fourier transform it holds for each $r \in {\mathbb N}_0$ that
$$
f^{(r)}(t) = \int_{-\delta}^{\delta} {\hat f}(v)\,(2 \pi {\mathrm i}v)^r\,{\mathrm e}^{2\pi {\mathrm i}v t}\,{\mathrm d}v\,,
$$
i.\,e., $f^{(r)} \in C_0(\mathbb R)$, because $(2 \pi {\mathrm i}\,\cdot)^r\,{\hat f}\in L^1([-\delta,\, \delta])$.
In the following we will always select the smooth representation of an equivalence class in ${\mathcal B}_{\delta}(\mathbb R)$.

In this paper we consider bandlimited functions $f \in {\mathcal B}_{\delta}(\mathbb R)$ with $\delta \in (0,\,N/2)$, where $N\in \mathbb N$ is fixed.
For \mbox{$L\coloneqq N(1+\lambda)$} with \mbox{$\lambda\geq 0$}, and any $m \in \mathbb N \setminus \{1\}$ with $2m \ll L$, we introduce the set $\Phi_{m,L}$ of all
window functions $\varphi: \mathbb R \to [0,\,1]$ with the following properties:
\begin{itemize}[leftmargin=*,nosep]
	\item[$\quad\bullet$]
	The window function \mbox{$\varphi\in L^2(\mathbb R)$} is even, positive on \mbox{$(-m/L,\,m/L)$} and continuous on \mbox{${\mathbb R}\setminus \{- m/L,\, m/L\}$}.
	\item[$\bullet$]
	The restricted window function $\varphi |_{[0,\,\infty)}$ is monotonously non-increasing with $\varphi(0) = 1\,.$
	\item[$\bullet$]
	The Fourier transform
	$$
	{\hat \varphi}(v) \coloneqq \int_{\mathbb R} \varphi(x)\,{\mathrm e}^{- 2 \pi {\mathrm i}v x}\,{\mathrm d}x = 2 \int_0^{\infty} \varphi(x)\,\cos( 2 \pi\,v x)\,{\mathrm d}x
	$$
	is explicitly known.
\end{itemize}
Examples of such window functions are the \emph{rectangular window function}
\begin{equation}
	\label{eq:varphirect}
	\varphi_{\mathrm{rect}} (x) \coloneqq {\mathbf 1}_{[-m/L,\,m/L]}(x)\,,\quad x \in \mathbb R\,,
\end{equation}
where ${\mathbf 1}_{[-m/L,\,m/L]}$ is the indicator function of the interval $[-m/L,\,m/L]$,
the \emph{Gaussian window function}
\begin{equation}
	\label{eq:varphiGauss}
	\varphi_{\mathrm{Gauss}}(x) \coloneqq {\mathrm e}^{- x^2/(2 \sigma^2)}\,, \quad x \in \mathbb R\,,
\end{equation}
with some $\sigma >0$, the modified $\mathrm B$-\emph{spline window function}
\begin{align}
	\label{eq:varphiB}
	\varphi_{\mathrm{B}}(x) \coloneqq \frac{1}{M_{2s}(0)}\, M_{2s}\bigg(\frac{Lxs}{m}\bigg)\,,\quad x \in \mathbb R\,,
\end{align}
where $M_{2s}$ is the centered cardinal $\mathrm B$--spline of even order $2s$, and the $\sinh$-\emph{type window function}
\begin{align}
	\label{eq:varphisinh}
	\varphi_{\sinh}(x) \coloneqq
	\begin{cases}
		\frac{1}{\sinh \beta}\, \sinh\big(\beta\,\sqrt{1-(L x/m)^2}\big) &\colon \quad x \in [-m/L,\,m/L] \,, \\
		0 &\colon \quad x \in \mathbb R\setminus [-m/L,\,m/L] \,,
	\end{cases}
\end{align}
with certain $\beta > 0$.
All these window functions are well-studied in the context of the nonequispaced fast Fourier transform (NFFT), see e.\,g. \cite{PT21a} and references therein.

Now let $\varphi \in \Phi_{m,L}$ be a given window function.
We define the \emph{truncated window function}
\begin{equation}
	\label{eq:varphim}
	\varphi_m (x) \coloneqq \varphi(x) \,{\mathbf 1}_{[-m/L,\,m/L]}(x)\,, \quad x \in \mathbb R\,,
\end{equation}
and study the \emph{regularized Shannon sampling formula with localized sampling}
\begin{equation}
	\label{eq:Rmf(t)}
	(R_{\varphi,m}f)(t) \coloneqq \sum_{\ell \in \mathbb Z} f\big(\tfrac{\ell}{L}\big)\,\mathrm{sinc}\big(L \pi\,\big(t-\tfrac{\ell}{L}\big)\big)\,\varphi_m\big(t - \tfrac{\ell}{L}\big)\,,
	\quad t \in \mathbb R\,,
\end{equation}
to rapidly reconstruct the values $f(t)$ for $t \in \mathbb R$ from given sampling data $f\big(\frac{\ell}{L}\big)$, $\ell \in \mathbb Z$, with high accuracy.
\medskip

It is known that $\{\mathrm{sinc}\big(L \pi\,\big(\cdot -\frac{\ell}{L}\big)\big) : \, \ell \in \mathbb Z\}$ forms an orthogonal system in $L^2(\mathbb R)$, since by the shifting property
the Fourier transform of $\mathrm{sinc}\big(L \pi\,\big(\cdot -\frac{\ell}{L}\big)\big)$ is equal to
$$
\frac{1}{L}\,{\mathrm e}^{- 2 \pi {\mathrm i}\ell v/L}\,{\mathbf 1}_{[-L/2,\,L/2]}(v)\,,\quad  v \in \mathbb R\,,
$$
and by the Parseval identity it holds for all $\ell$, $k \in \mathbb Z$ that
\begin{equation}
\label{eq:orthsystem}
	\big\langle \mathrm{sinc}\big(L \pi\,\big( \cdot -\tfrac{\ell}{L}\big)\big),\,\mathrm{sinc}\big(L \pi\,\big( \cdot -\tfrac{k}{L}\big)\big)\big\rangle_{L^2(\mathbb R)}= \frac{1}{L^2}\, \int_{-L/2}^{L/2} {\mathrm e}^{2 \pi {\mathrm i}(k-\ell)\,v/L}\,{\mathrm d}v = \frac{1}{L}\, \delta_{\ell,k}
\end{equation}
with the Kronecker symbol $\delta_{\ell,k}$.
Moreover, it follows directly from the Whittaker--Kotelnikov--Shannon sampling theorem (see \eqref{eq:WKSseries}) that the system
\mbox{$\{\mathrm{sinc}\big(L \pi\,\big(\cdot -\frac{\ell}{L}\big)\big) : \, \ell \in \mathbb Z\}$}
forms an orthogonal basis of ${\mathcal B}_{N/2}(\mathbb R)$.
From \eqref{eq:WKSseries} and \eqref{eq:orthsystem} it follows that
for any \mbox{$f \in {\mathcal B}_{\delta}(\mathbb R) \subset {\mathcal B}_{L/2}(\mathbb R)$} with $\delta \in (0,\, N/2]$ and $L\geq N$ it holds
\begin{equation}
	\label{eq:Parseval}
	\| f \|_{L^2(\mathbb R)}^2 = \tfrac{1}{L}\, \sum_{\ell \in \mathbb Z} \big|f\big(\tfrac{\ell}{L}\big)\big|^2\,.
\end{equation}

Firstly, we consider the regularized Shannon sampling formula \eqref{eq:Rmf(t)} with the simple rectangular window function $\varphi = \varphi_{\mathrm{rect}}$, see~\eqref{eq:varphirect}, i.\,e., for some $m \in \mathbb N \setminus \{1\}$ we form the \emph{rectangular regularized Shannon sampling formula with localized sampling}
\begin{equation}
\label{eq:Rmrectf}
(R_{{\mathrm{rect}},m}f)(t) \coloneqq \sum_{\ell \in \mathbb Z} f\big(\tfrac{\ell}{L}\big) \, \mathrm{sinc}\big(L\pi\,\big(t-\tfrac{\ell}{L}\big)\big)\,{\mathbf 1}_{[-m/L,\,m/L]}\big(t - \tfrac{\ell}{L}\big)\,, \quad t \in \mathbb R\,.
\end{equation}
Obviously, the rectangular regularized Shannon sampling formula \eqref{eq:Rmrectf} interpolates $f$ on the grid $\frac{1}{L}\,\mathbb Z$, i.\,e.,
for all $m\in \mathbb N\setminus \{1\}$, the interpolation property
\begin{equation}
\label{eq:(Rmf)(j/N)}
f\big(\tfrac{\ell}{L}\big) = (R_{{\mathrm{rect}},m} f)\big(\tfrac{\ell}{L}\big)\,, \quad \ell \in \mathbb Z\,,
\end{equation}
is fulfilled since $\sinc(\pi k)=0$ for $k\in\Z\setminus\{0\}$.
Due to the definition of the indicator function ${\mathbf 1}_{[-m/L,\,m/L]}$, for $t \in \big(0,\,\frac{1}{L}\big)$ the rectangular regularized Shannon sampling formula reads as
$$
(R_{{\mathrm{rect}},m} f)(t) = \sum_{\ell\in\mathcal J_m} f\big(\tfrac{\ell}{L}\big) \, \mathrm{sinc}\big(L\pi\,\big(t - \tfrac{\ell}{L}\big)\big)
$$
with the index set \mbox{$\mathcal J_m \coloneqq \{-m+1,\,-m+2,\,\ldots,\,m\}$}.
Indeed, on any interval $\big(\frac{k}{L},\,\frac{k+1}{L}\big)$ with \mbox{$k\in \mathbb Z$} the rectangular regularized Shannon sampling formula reads as
\begin{equation}
\label{eq:Smfinterval}
(R_{{\mathrm{rect}},m} f)\big(t+\tfrac{k}{L}\big) = \sum_{\ell\in\mathcal J_m} f\big(\tfrac{\ell+k}{L}\big) \, \mathrm{sinc}\big(L\pi\,\big(t - \tfrac{\ell}{L}\big)\big)\,, \quad t \in \big(0,\,\tfrac{1}{L}\big) \,.
\end{equation}

However, since the $\sinc$ function decays slowly at infinity, \eqref{eq:Rmrectf} is not a good approximation to $f$ on~$\mathbb R$.
As a consequence of a result in \cite{MXZ09}, it can be seen that the convergence rate of the sequence $\big(f - R_{{\mathrm{rect}},m} f \big)_{m=1}^{\infty}$ is only ${\mathcal{O}}(m^{-1/2})$ for sufficiently large $m$.

\begin{Lemma}
\label{Lemma:f-Rmrectf}
	Let $f \in {\mathcal{B}}_{N/2}(\mathbb R)$ with fixed $N \in \mathbb N$, $L\coloneqq N(1+\lambda)$ with \mbox{$\lambda\geq 0$} and $m \in \mathbb N \setminus\{1\}$ be given. Then it holds
	$$
	\| f - R_{{\mathrm{rect}},m} f\|_{C_0(\mathbb R)} \le \frac{\sqrt{L}}{\pi}\,\sqrt{\frac{2}{m} + \frac{1}{m^2}}\,\|f\|_{L^2(\mathbb R)}\,.
	$$
\end{Lemma}

\emph{Proof.} Since $R_{{\mathrm{rect}},m} f$ possesses similar representations \eqref{eq:Smfinterval} on each interval  $\big(\frac{k}{L},\,\frac{k+1}{L}\big)$, $k\in \mathbb Z$, we consider $f(t) - (R_{m,{\mathrm{rect}}} f)(t)$ only for $t \in \big[0,\,\frac{1}{L}\big]$ and show that
\begin{equation}
\label{eq:|f-Smf|interval}
\max_{t\in [0,\,1/L]} |f(t) - (R_{{\mathrm{rect}},m} f)(t)| \le \frac{\sqrt{L}}{\pi}\,\sqrt{\frac{2}{m} + \frac{1}{m^2}}\,\|f\|_{L^2(\mathbb R)}\,.
\end{equation}
The Whittaker--Kotelnikov--Shannon sampling theorem (see \eqref{eq:WKSseries}) implies that
$$
f(t) - (R_{{\mathrm{rect}},m} f)(t) = \sum_{\ell \in \Z\setminus\mathcal J_m} f\big(\tfrac{\ell}{L}\big) \,{\mathrm{sinc}}\big(L \pi\,\big(t - \tfrac{\ell}{L}\big)\big)\,.
$$
Then by the Cauchy--Schwarz inequality and \eqref{eq:Parseval} it follows that
\begin{eqnarray}
\label{eq:|f-Snf|le}
	\big| f(t) - (R_{{\mathrm{rect}},m} f)(t) \big|
	&\le& \Big(\sum_{\ell \in \Z\setminus\mathcal J_m} \big|f\big(\tfrac{\ell}{L}\big)\big|^2\Big)^{1/2}\,\Big(\sum_{\ell \in \Z\setminus\mathcal J_m}
	\big[{\mathrm{sinc}}\big(L \pi \,\big(t-\tfrac{\ell}{L}\big)\big)\big]^2\Big)^{1/2} \nonumber \\
	&\le& \sqrt{L}\,\| f\|_{L^2(\mathbb R)}\,\sqrt{h_m(t)}
\end{eqnarray}
with the auxiliary function
$$
h_m(t) \coloneqq \sum_{\ell \in \Z\setminus\mathcal J_m} \big[{\mathrm{sinc}}\big(L\pi \,\big(t-\tfrac{\ell}{L}\big)\big)\big]^2
= \frac{(\sin (L\pi t))^2}{\pi^2}\, \sum_{\ell \in {\mathbb Z}\setminus \mathcal J_m} \frac{1}{(L t-\ell)^2}
\ge 0\,, \quad t \in \big[0,\, \tfrac{1}{L}\big]\,.
$$
By the integral test for convergence of series we estimate the function
\begin{eqnarray}
	\label{eq:boundh(t)}
	h_m(t)
	\le \frac{1}{\pi^2}\, \sum_{k\in {\mathbb Z}\setminus \mathcal J_m} \frac{1}{k^2}
	\le \frac{1}{\pi^2}\,\left( \frac{1}{m^2} + 2 \, \int_m^{\infty} \frac{1}{x^2} \, \mathrm{d}x\right) =
	\frac{1}{\pi^2}\,\left( \frac{1}{m^2} + \frac{2}{m}\right)\,
\end{eqnarray}
for $t \in \big(0,\, \frac{1}{L}\big)$.
Then \eqref{eq:|f-Snf|le}, \eqref{eq:boundh(t)} combined with \eqref{eq:(Rmf)(j/N)} imply \eqref{eq:|f-Smf|interval}.

By the same technique, the above estimate of the approximation error
can be shown
on each interval  $\big(\frac{k}{L},\,\frac{k+1}{L}\big)$, $k\in \mathbb Z$.
This completes the proof. \qedsymbol
\medskip

In view of the slow convergence of the sequence $\big(R_{{\mathrm{rect}},m} f(t)\big)_{m=1}^{\infty}$ it has been proposed
to modify the rectangular regularized Shannon sampling sum \eqref{eq:Rmrectf} by multiplying the $\sinc$ function with a more convenient
window function $\varphi \in \Phi_{m,L}$ (see \cite{Q03} and \cite{LZ16}).
For any $m \in \mathbb N \setminus \{1\}$ the {regularized Shannon sampling formula with localized sampling} is given by
\begin{equation}
\label{eq:Rmf}
	(R_{\varphi,m}f)(t) = \sum_{\ell\in \mathbb Z}  f\big(\tfrac{\ell}{L}\big)\,\mathrm{sinc}\big(L \pi\,\big(t-\tfrac{\ell}{L}\big)\big) \, \varphi_m\big(t - \tfrac{\ell}{L}\big)\,,
	\quad t \in \mathbb R\,,
\end{equation}
with the {truncated window function} \eqref{eq:varphim}.
Note that it holds the interpolation property
\begin{equation}
\label{eq:(Rmf)Z}
	f\big(\tfrac{\ell}{L}\big) = (R_{\varphi, m}f)\big(\tfrac{\ell}{L}\big)\,, \quad \ell \in \mathbb Z\,.
\end{equation}
Especially for $t \in \big(0,\,\frac{1}{L}\big)$, we obtain the regularized Shannon sampling formula
$$
(R_{\varphi,m} f)(t)
=
\sum_{\ell\in {\mathcal J_m}} f\big(\tfrac{\ell}{L}\big)\,\mathrm{sinc}\big(L \pi\,\big(t-\tfrac{\ell}{L}\big)\big) \, \varphi_m\big(t - \tfrac{\ell}{L}\big)
=
\sum_{\ell\in {\mathcal J_m}} f\big(\tfrac{\ell}{L}\big)\, \psi\big(t - \tfrac{\ell}{L}\big) \,,
$$
where
\begin{equation}
\psi(x) \coloneqq \sinc(L\pi x) \,\varphi(x)  \label{eq:psi}
\end{equation}
is the \textit{regularized $\sinc$ function}.
For the reconstruction of $f$ on any interval $\big(\frac{k}{L},\,\frac{k+1}{L}\big)$ with $k \in \mathbb Z$, we use
\begin{equation}
\label{eq:Rmfintervalk}
(R_{\varphi,m} f)\big(t + \tfrac{k}{L}\big) = \sum_{\ell\in {\mathcal J_m}} f\big(\tfrac{\ell+k}{L}\big)\, \psi\big(t - \tfrac{\ell}{L}\big)\,,
\quad t \in \big(0,\,\tfrac{1}{L}\big)\,,
\end{equation}
i.\,e., we reconstruct $f$ by $R_{\varphi,m}f$ separately for each open interval $\big(\frac{k}{L},\,\frac{k+1}{L}\big)$, $k\in \mathbb Z$.
Now we estimate the uniform approximation error
\begin{align}
\label{eq:Error_RWKS}
	\| f - R_{\varphi,m} f \|_{C_0(\mathbb R)} \coloneqq \max_{t\in \mathbb R} \big| f(t) - (R_{\varphi,m}f)(t) \big|
\end{align}
of the regularized Shannon sampling formula.

\begin{Theorem}
\label{Theorem:errorRmf}
Let $f\in {\mathcal B}_{\delta}(\mathbb R)$ with $\delta = \tau N$, $\tau\in (0,\,1/2)$, $N \in \mathbb N$, \mbox{$L= N(1+\lambda)$} with \mbox{$\lambda\geq 0$} and
\mbox{$m \in {\mathbb N}\setminus \{1\}$}.
Further let $\varphi \in \Phi_{m,L}$ with the truncated window function \eqref{eq:varphim} be given. \\
Then the regularized Shannon sampling formula \eqref{eq:Rmf} with localized sampling satisfies
\begin{align}
\label{eq:error_const}
	\| f - R_{\varphi,m}f \|_{C_0(\mathbb R)} \le \big( E_1(m,\delta,L) +  E_2(m,\delta,L) \big) \,\|f\|_{L^2(\mathbb R)} \,,
\end{align}
where the corresponding error constants are defined by
\begin{align}
E_1(m,\delta,L) &\coloneqq \sqrt{2\delta} \, \max_{v \in [-\delta,\delta]} \left| 1 - \int_{v - \frac L2}^{v + \frac L2} {\hat\varphi}(u)\,{\mathrm d}u \,\right|\,, \label{eq:E1} \\
E_2(m,\delta,L) &\coloneqq \frac{\sqrt{2L}}{\pi m}\, \left(\varphi^2\big(\tfrac{m}{L}\big) + L\int_{\frac{m}{L}}^{\infty} \varphi^2(t)\,{\mathrm d}t\right)^{1/2} \,. \label{eq:E2}
\end{align}
\end{Theorem}

\emph{Proof}.
By \eqref{eq:Rmfintervalk} we split the appro\-xi\-mation error
$$
f\big(t + \tfrac{k}{L}\big) - (R_{\varphi,m} f)\big(t + \tfrac{k}{L}\big) = e_1\big(t + \tfrac{k}{L}\big) + e_{2,k}(t)\,, \quad t \in \big[0,\, \tfrac{1}{L}\big]\,,
$$
on each interval $\big[\frac{k}{L},\,\frac{k+1}{L}\big]$ with $k \in \mathbb Z$ into the regularization error
\begin{align}
\label{eq:e1}
	e_1\big(t + \tfrac{k}{L}\big) &\coloneqq f\big(t + \tfrac{k}{L}\big) - \sum_{\ell\in \mathbb Z} f\big(\tfrac{\ell + k}{L}\big)\,\psi\big(t - \tfrac{\ell}{L}\big) \,, \quad t \in \big[0,\, \tfrac{1}{L}\big] \,,
\end{align}
where $\psi$ denotes the regularized $\mathrm{sinc}$ function \eqref{eq:psi},
and the truncation error
\begin{align}
\label{eq:e2}
	e_{2,k}(t)
	&\coloneqq
	\sum_{\ell \in \mathbb Z} f\big(\tfrac{\ell + k}{L}\big)\, \psi\big(t - \tfrac{\ell}{L}\big) - (R_{\varphi,m} f)(t) \,, \quad t \in \big[0,\,\tfrac{1}{L}\big] \,.
\end{align}

\vspace*{-1ex}
Initially, we only consider the error on the interval $\big[0,\,\frac{1}{L}\big]$.
We start with the regularization error \eqref{eq:e1}.
By Lemma \ref{Lemma:FTconvolution}, the Fourier transform of $\psi$ reads as
$$
{\hat\psi}(v) = \frac 1L \int_\R {\mathbf 1}_{[-L/2,\,L/2]}(v-u) \, {\hat\varphi}(u) \,{\mathrm d}u  = \frac 1L \int_{v-L/2}^{v+L/2} {\hat\varphi}(u) \,{\mathrm d}u\,.
$$
Hence, using the shifting property of the Fourier transform, the Fourier transform of $\psi\big(\cdot \, - \frac{\ell}{L}\big)$ reads as
\vspace*{-2.5ex}
\begin{eqnarray*}
	\frac{1}{L}\,{\mathrm e}^{-{2 \pi \mathrm i} v \ell/L}\,\int_{v - L/2}^{v + L/2} \hat\varphi(u)\,{\mathrm d}u\,.
\end{eqnarray*}
Therefore, the Fourier transform of the regularization error $e_1$ has the form
\begin{align}
\label{eq:e1_FT}
	{\hat e}_1(v) = {\hat f}(v) - \left(\sum_{\ell\in \mathbb Z} f\big(\tfrac{\ell}{L}\big)\,\frac{1}{L}\,{\mathrm e}^{-2 \pi {\mathrm i} v \ell/L}\right)\,\int_{v - L/2}^{v + L/2}
	{\hat\varphi}(u)\,{\mathrm d}u\,.
\end{align}
By the assumption $f \in {\mathcal B}_{\delta}(\mathbb R)$ with $\delta \in (0,\,N/2)$ and \mbox{$L\geq N$}, it holds
$
\mathrm{supp}\,\hat f \subseteq [- \delta,\, \delta] \subset [- L/2,\,L/2 ]
$
and hence the restricted function
\mbox{$\hat f|_{[-L/2,\,L/2]}$} belongs to \mbox{$L^2([-L/2,\,L/2])$}.
Thus, this function possesses the $L$-periodic Fourier expansion
$$
{\hat f}(v) = \sum_{k \in \mathbb Z} c_k(\hat f)\, {\mathrm e}^{-2 \pi{\mathrm i} k v/L}\,, \quad v \in [-L/2, \,L/2]\,,
$$
with the Fourier coefficients
$$
c_k(\hat f) = \frac{1}{L}\,\int_{-L/2}^{L/2} {\hat f}(u)\,{\mathrm e}^{2\pi {\mathrm i} k u/L}\, {\mathrm d}u = \frac{1}{L}\,f\big(\tfrac{k}{L}\big) \,, \quad k \in \Z \,.
$$
In other words, ${\hat f}$ can be represented as
\begin{align}
\label{eq:f_FT}
{\hat f}(v) = {\hat f}(v)\,{\mathbf 1}_{[-\delta,\,\delta]}(v) = \Big(\sum_{k \in \mathbb Z} \tfrac{1}{L}\,f\big(\tfrac{k}{L}\big)\,{\mathrm e}^{-2 \pi {\mathrm i}k v/L}\Big)
\,{\mathbf 1}_{[-\delta,\,\delta]}(v) \,, \quad v \in \R \,.
\end{align}
Introducing the auxiliary function
\begin{align}
\label{eq:eta}
	\eta(v) \coloneqq {\mathbf 1}_{[-\delta,\,\delta]}(v) - \int_{v - L/2}^{v + L/2} {\hat\varphi}(u)\,{\mathrm d}u \,, \quad v \in \R \,,
\end{align}
we see by inserting \eqref{eq:f_FT} into \eqref{eq:e1_FT} that
\mbox{${\hat e}_1(v) = {\hat f}(v)\, \eta (v)$}
and thereby
\mbox{$|{\hat e}_1(v)| \le |{\hat f}(v)|\,|\eta(v)|$}.
Thus, by inverse Fourier transform \eqref{eq:inverse_Fouriertrafo} we get
\begin{align*}
	|e_1(t)|
	\leq \int_{\mathbb R} |{\hat e}_1(v)|\,{\mathrm d}v
	\leq \int_{-\delta}^{\delta} |{\hat f}(v)|\,|\eta(v)|\,{\mathrm d}v
	\leq \max_{v \in [-\delta,\delta]} \,|\eta(v)| \,\int_{-\delta}^{\delta} |{\hat f}(v)|\,{\mathrm d}v \,.
\end{align*}
Using Cauchy--Schwarz inequality and Parseval identity, we see that
\begin{eqnarray*}
	\int_{-\delta}^{\delta} |{\hat f}(v)|\,{\mathrm d}v &\le& \left(\int_{-\delta}^{\delta} 1^2\,{\mathrm d}v\right)^{1/2} \left(\int_{-\delta}^{\delta} |{\hat f}(v)|^2\,{\mathrm d}v\right)^{1/2}
	= \sqrt{2 \delta}\,\|{\hat f}\|_{L^2(\mathbb R)} = \sqrt{2 \delta}\,\|f\|_{L^2(\mathbb R)}.
\end{eqnarray*}
In summary, using the error constant \eqref{eq:E1}, this yields
$
\|e_1 \|_{C_0(\mathbb R)} \le E_1(m,\delta,L)\,\|f\|_{L^2(\mathbb R)} .
$

Now we estimate the truncation error. By \eqref{eq:e2} it holds for $t \in \big(0,\,\frac{1}{L}\big)$ that
\begin{align*}
	e_{2,0}(t)
	&=
	\sum_{\ell \in \mathbb Z} f\big(\tfrac{\ell}{L}\big)\, \psi\big(t - \tfrac{\ell}{L}\big) \,\big[1- {\mathbf 1}_{[-m/L,m/L]}\big(t - \tfrac{\ell}{L}\big)\big]
	= \sum_{\ell \in \mathbb Z \setminus \mathcal J_m} f\big(\tfrac{\ell}{L}\big)\, \psi\big(t - \tfrac{\ell}{L}\big) \,.
\end{align*}
Using \eqref{eq:psi}  and the non-negativity of $\varphi$, we receive
$$
|e_{2,0}(t)| \le \sum_{\ell \in \mathbb Z \setminus {\mathcal J_m}} \big|f\big(\tfrac{\ell}{L}\big)\big|\, \big|\mathrm {sinc}\big(L \pi \,\big(t - \tfrac{\ell}{L}\big)\big)\big| \,\varphi\big(t - \tfrac{\ell}{L}\big)\,.
$$
For $t \in \big(0,\, \frac{1}{L}\big)$ and $\ell \in {\mathbb Z} \setminus \mathcal J_m$ we obtain
$$
\big|\mathrm {sinc}\big(L \pi \,\big(t - \tfrac{\ell}{L}\big)\big)\big| \le \frac{1}{\pi \,|L t - \ell|} \le \frac{1}{\pi m}
$$
and hence
$$
|e_{2,0}(t)| \le  \frac{1}{\pi m}\, \sum_{\ell \in \mathbb Z \setminus \mathcal J_m} \big|f\big(\tfrac{\ell}{L}\big)\big|\, \varphi\big(t - \tfrac{\ell}{L}\big)\,.
$$
Then the Cauchy--Schwarz inequality implies that
$$
|e_{2,0}(t)| \le  \frac{1}{\pi m}\, \bigg( \sum_{\ell \in \mathbb Z \setminus \mathcal J_m} \big|f\big(\tfrac{\ell}{L}\big)\big|^2 \bigg)^{1/2}\bigg(\sum_{\ell \in \mathbb Z \setminus \mathcal J_m} \varphi^2\big(t - \tfrac{\ell}{L}\big) \bigg)^{1/2}\,.
$$
By \eqref{eq:Parseval} it holds
$$
\bigg( \sum_{\ell \in \mathbb Z \setminus \mathcal J_m} \big|f\big(\tfrac{\ell}{L}\big)\big|^2 \bigg)^{1/2} \le \sqrt L\, \|f \|_{L^2(\mathbb R)}\,.
$$
Since $\varphi |_{[0,\,\infty)}$ is monotonously non-increasing by assumption $\varphi\in\Phi_{m,L}$, we can estimate the series for $t \in \big(0,\, \frac{1}{L}\big)$ by
\begin{align*}
	\sum_{\ell \in \mathbb Z \setminus \mathcal J_m} \varphi^2\big(t - \tfrac{\ell}{L}\big)
	&= \Bigg(\sum_{\ell=-\infty}^{-m} + \sum_{\ell=m+1}^{\infty}\Bigg)\, \varphi^2\big(t - \tfrac{\ell}{L}\big)
	= \sum_{\ell=m}^{\infty} \,\varphi^2\big(t + \tfrac{\ell}{L}\big) + \sum_{\ell=m+1}^{\infty} \,\varphi^2\big(t - \tfrac{\ell}{L}\big) \\
	&\leq \sum_{\ell=m}^{\infty} \,\varphi^2\big(\tfrac{\ell}{L}\big) + \sum_{\ell=m+1}^{\infty} \,\varphi^2\big(\tfrac 1L - \tfrac{\ell}{L}\big)
	= 2 \sum_{\ell=m}^{\infty} \,\varphi^2\big(\tfrac{\ell}{L}\big) \,.
\end{align*}
Using the integral test for convergence of series, we obtain that
\begin{align*}
	\sum_{\ell=m}^{\infty} \,\varphi^2\big(\tfrac{\ell}{L}\big)
	= \varphi^2\big(\tfrac{m}{L}\big) + \sum_{\ell=m+1}^{\infty} \,\varphi^2\big(\tfrac{\ell}{L}\big)
	< \varphi^2\big(\tfrac{m}{L}\big) + \int_{m}^{\infty} \varphi^2\big(\tfrac{t}{L}\big)\,{\mathrm d}t
	= \varphi^2\big(\tfrac{m}{L}\big) + L \, \int_{m/L}^{\infty} \varphi^2(t)\,{\mathrm d}t \,.
\end{align*}
By \eqref{eq:(Rmf)Z} it holds $e_{2,0}(0) = e_{2,0}\big(\frac{1}{L}\big) = 0$.
Hence, we obtain by \eqref{eq:E2} that
\begin{equation*}
	\max_{t \in [0,1/L]} |e_{2,0}(t)| \le \frac{\sqrt{L}}{\pi m}\, \|f \|_{L^2(\mathbb R)} \left(2\varphi^2\big(\tfrac{m}{L}\big) + 2L\int_{m/L}^{\infty} \varphi^2(t)\,{\mathrm d}t\right)^{1/2}
	= E_2(m,\delta,L) \, \|f \|_{L^2(\mathbb R)} \,.
\end{equation*}

By the same technique, this error estimate can be shown for each interval $\big[\frac{k}{L},\,\frac{k+1}{L}\big]$ with $k \in \mathbb Z$.
This completes the proof.
\qedsymbol
\medskip

\begin{Remark}
\label{Remark:simpleerror}
Theorem \ref{Theorem:errorRmf} can be simplified, if the window function $\varphi \in \Phi_{m,L}$ is continuous on $\mathbb R$ and vanishes on $\mathbb R \setminus \big[- \frac{m}{L},\, \frac{m}{L}\big]$.
Then the truncation errors $e_{2,k}(t)$ vanish for all $t \in \big(0,\,\frac{1}{L}\big)$ and $k \in \mathbb Z$, such that $E_2(m,\delta,L) = 0$.
Thereby, we obtain the simple error estimate
$$
\| f - R_{\varphi,m}f \|_{C_0(\mathbb R)} \le  E_1(m,\delta,L)\,\|f\|_{L^2(\mathbb R)}\,.
$$
We remark that this is the case for the $\mathrm B$--spline~\eqref{eq:varphiB} as well as the $\sinh$-type window function~\eqref{eq:varphisinh}, but not for the Gaussian window function~\eqref{eq:varphiGauss} since $\varphi_{\mathrm{Gauss}}$ does not vanish on \mbox{$\mathbb R \setminus \big[- \frac{m}{L},\, \frac{m}{L}\big]$}.
Also the rectangular window function~\eqref{eq:varphirect} does not fit into this setting since $\varphi_{\mathrm{rect}}$ is not continuous on~$\R$.
\end{Remark}

If the samples $f\big(\frac{\ell}{L}\big)$, $\ell \in \mathbb Z$, of a bandlimited function $f\in {\mathcal B}_{\delta}(\mathbb R)$ are not known exactly, i.\,e.,
only erroneous samples ${\tilde f}_{\ell} \coloneqq f\big(\frac{\ell}{L}\big) + \varepsilon_{\ell}$ with $|\varepsilon_{\ell}| \le \varepsilon$, $\ell\in\Z$, with $\varepsilon>0$ are known,
the corresponding Shannon sampling series may differ appreciably from $f$ (see \cite{DDeV03}).
Here we denote the regularized Shannon sampling formula with erroneous samples ${\tilde f}_{\ell}$ by
\begin{align}
\label{eq:Rmf_erroneous}
	(R_{\varphi,m}{\tilde f})(t) = \sum_{\ell\in \mathbb Z}  \,\tilde f_\ell\ \mathrm{sinc}\big(L \pi\,\big(t-\tfrac{\ell}{L}\big)\big) \, \varphi_m\big(t - \tfrac{\ell}{L}\big)\,,
	\quad t \in \mathbb R\,.
\end{align}
In contrast to the classical Shannon sampling series, the regularized Shannon sampling formula is numerically robust, i.\,e., the uniform perturbation error
\begin{align}
\label{eq:perturb_error}
	\| R_{\varphi,m}{\tilde f} - R_{\varphi,m}f \|_{C_0(\mathbb R)}
	\coloneqq \max_{t\in \mathbb R} \big| (R_{\varphi,m}{\tilde f})(t) - (R_{\varphi,m}f)(t) \big|
\end{align}
is small, as shown in the following.

\begin{Theorem}
\label{Theorem:robustness}
	Let $f\in {\mathcal B}_{\delta}(\mathbb R)$ with $\delta = \tau N$, $\tau\in (0,\,1/2)$, $N \in \mathbb N$, \mbox{$L= N(1+\lambda)$} with \mbox{$\lambda\geq 0$} and \mbox{$m \in {\mathbb N}\setminus \{1\}$} be given.
	Further let $\varphi \in \Phi_{m,L}$ with the truncated window function \eqref{eq:varphim} as well as ${\tilde f}_{\ell} = f(\ell/L) + \varepsilon_{\ell}$,
	where $|\varepsilon_{\ell}| \le \varepsilon$ for all $\ell\in\Z$, with $\varepsilon>0$ be given.\\
	Then the regularized Shannon sampling sum \eqref{eq:Rmf} with localized sampling satisfies
	\begin{eqnarray}
		\| R_{\varphi,m}{\tilde f} - R_{\varphi,m}f \|_{C_0(\mathbb R)} &\le& \varepsilon \, \big( 2+L\,\hat\varphi(0) \big) \,, \label{eq:result_robustness} \\
		\| f - R_{\varphi,m}{\tilde f} \|_{C_0(\mathbb R)} &\le& \| f - R_{\varphi,m}{f}\|_{C_0(\mathbb R)} + \varepsilon \, \big( 2+L\,\hat\varphi(0) \big) \,. \label{eq:robust2}
	\end{eqnarray}
\end{Theorem}

\emph{Proof}.
On each open interval $\big(\frac{k}{L},\,\frac{k+1}{L}\big)$ with $k \in \mathbb Z$ we denote the error by \eqref{eq:Rmfintervalk} as
\begin{align*}
	\tilde e_{k}(t)
	\coloneqq (R_{\varphi,m} {\tilde f})\big(t + \tfrac{k}{L}\big) - (R_{\varphi,m} f)\big(t + \tfrac{k}{L}\big)
	=\sum_{\ell \in \mathcal J_m} \varepsilon_{\ell + k}\,\psi\big(t - \tfrac{\ell}{L}\big)
	\,, \quad t \in \big(0,\, \tfrac{1}{L}\big)\,.
\end{align*}
Initially, we consider the interval $\big[0,\,\frac{1}{L}\big]$.
Using \eqref{eq:psi}, the non-negativity of $\varphi$ and $|\varepsilon_{\ell}| \le \varepsilon$, we receive
\begin{align*}
	|\tilde e_0(t)|
	&\leq
	\sum_{\ell \in \mathcal J_m} \left|\varepsilon_{\ell}\right|\, \big|\mathrm {sinc}\big(L \pi \,\big(t - \tfrac{\ell}{L}\big)\big)\big| \,\varphi\big(t - \tfrac{\ell}{L}\big)
	\leq
	\varepsilon \sum_{\ell \in \mathcal J_m} \,\varphi\big(t - \tfrac{\ell}{L}\big) \,.
\end{align*}
Since $\varphi |_{[0,\,\infty)}$ is monotonously non-increasing by assumption $\varphi\in\Phi_{m,L}$, we can estimate the sum for $t \in \big(0,\, \frac{1}{L}\big)$ by
\begin{align*}
\sum_{\ell \in \mathcal J_m} \varphi\big(t - \tfrac{\ell}{L}\big)
&= \Big(\sum_{\ell=-m+1}^{0} + \sum_{\ell=1}^{m}\Big)\, \varphi\big(t - \tfrac{\ell}{L}\big)
= \sum_{\ell=0}^{m-1} \,\varphi\big(t + \tfrac{\ell}{L}\big) + \sum_{\ell=1}^{m} \,\varphi\big(t - \tfrac{\ell}{L}\big) \\
&\leq \sum_{\ell=0}^{m-1} \,\varphi\big(\tfrac{\ell}{L}\big) + \sum_{\ell=1}^{m} \,\varphi\big(\tfrac 1L - \tfrac{\ell}{L}\big)
= 2 \sum_{\ell=0}^{m-1} \,\varphi\big(\tfrac{\ell}{L}\big) \,.
\end{align*}
Using the integral test for convergence of series, we obtain that
\begin{align*}
	\sum_{\ell=0}^{m-1} \,\varphi\big(\tfrac{\ell}{L}\big)
	< \varphi(0) + \int_{0}^{m-1} \varphi\big(\tfrac{t}{L}\big)\,{\mathrm d}t
	= \varphi(0) + L \, \int_{0}^{(m-1)/L} \varphi(t)\,{\mathrm d}t \,.
\end{align*}
By the definition of the Fourier transform \eqref{eq:Fouriertrafo} it holds for $\varphi \in \Phi_{m,L}$ that
\begin{align*}
	{\hat \varphi}(0)
	=
	\int_{\mathbb R} \varphi(t)\, {\mathrm d}t
	\geq
	\int_{-m/L}^{m/L} \varphi(t)\, {\mathrm d}t
	=
	2\int_{0}^{m/L} \varphi(t)\, {\mathrm d}t
	\geq
	2\int_{0}^{(m-1)/L} \varphi(t)\, {\mathrm d}t \,,
\end{align*}
and therefore
\begin{align*}
	|\tilde e_0(t)|
	&\leq
	2\,\varepsilon \sum_{\ell=0}^{m-1} \,\varphi\big(\tfrac{\ell}{L}\big)
	\leq
	2\,\varepsilon \big( \varphi(0) + \tfrac L2 \, {\hat \varphi}(0) \big)
	=
	\varepsilon \big( 2\,\varphi(0) + L \, {\hat \varphi}(0) \big)
	\,, \quad t \in \big(0,\, \tfrac{1}{L}\big) \,.
\end{align*}
Additionally, by the interpolation property \eqref{eq:(Rmf)Z} it holds
$|\tilde e_0(0)| = |\varepsilon_0| \leq \varepsilon$ as well as
$\left|\tilde e_0\big(\frac{1}{L}\big)\right| = |\varepsilon_1| \leq \varepsilon$.
By $\varphi \in \Phi_{m,L}$ we have $\varphi(0)=1$ and therefore we obtain that
\begin{equation*}
\max_{t \in [0,1/L]} |\tilde e_0(t)| \leq \varepsilon \big( 2 + L \, {\hat \varphi}(0) \big) \,.
\end{equation*}

By the same technique, this error estimate can be shown for each interval $\big[\frac{k}{L},\,\frac{k+1}{L}\big]$ with $k \in \mathbb Z$.
Then the triangle inequality yields \eqref{eq:robust2}, which completes the proof.
\qedsymbol
\medskip

Now it merely remains to estimate the error constants \mbox{$E_j(m,\delta,L)$}, \mbox{$j=1,2,$} for a certain window function, which shall be done for some selected ones in the following sections.

\section{Gaussian regularized Shannon sampling formula\label{sec:Gaussian}}

Firstly, we consider the Gaussian window function~\eqref{eq:varphiGauss}
with some $\sigma >0$
and show that in this case the uniform approximation error~\eqref{eq:Error_RWKS} for the regularized Shannon sampling formula~\eqref{eq:Rmf} decays exponentially with respect to $m$.

\begin{Theorem}
	\label{Theorem:GaussianWKS}
	Let $f\in {\mathcal B}_{\delta}(\mathbb R)$ with $\delta = \tau N$, $\tau \in (0,\,1/2)$, $N \in \mathbb N$, \mbox{$L= N(1+\lambda)$} with \mbox{$\lambda\geq 0$}, and \mbox{$m \in {\mathbb N}\setminus \{1\}$} be given.\\
	Then the regularized Shannon sampling formula \eqref{eq:Rmf} with the Gaussian window function~\eqref{eq:varphiGauss} and $\sigma = \sqrt{\frac{m}{\pi L\,(L - 2\delta)}}$ satisfies the error estimate
	\begin{align}
		\label{eq:error_result_Gaussian}
		\| f - R_{\mathrm{Gauss},m}f \|_{C_0(\mathbb R)} \le \frac{2\sqrt{\pi \delta L} + L(m+1)/\sqrt m}{\pi\,\sqrt{m\pi(L- 2 \delta)}}\,{\mathrm e}^{-\pi m (L/2 -
			\delta)/L}\,\|f\|_{L^2(\mathbb R)}\,.
	\end{align}
\end{Theorem}

\emph{Proof} (cf. \cite{Q03} and \cite{LZ16}). By Theorem \ref{Theorem:errorRmf} we have to compute the error constants \mbox{$E_j(m,\delta,L)$}, \mbox{$j=1,2,$} for the Gaussian window function \eqref{eq:varphiGauss}.
First we study the regularization error constant \eqref{eq:E1}.
The function \eqref{eq:varphiGauss} possesses the Fourier transform
\begin{align}
	\label{eq:Gaussian_FT}
	{\hat \varphi_{\mathrm{Gauss}}}(v) &= \sqrt{2\pi}\,\sigma\,{\mathrm e}^{- 2\pi^2 \sigma^2 v^2}\,, \quad v \in \mathbb R\,.
\end{align}
Thus, we obtain by substitution $w=\sqrt 2 \pi \sigma u$ that the auxiliary function \eqref{eq:eta} is given by
$$
\eta(v) = {\mathbf 1}_{[-\delta,\,\delta]}(v) - \frac{1}{\sqrt \pi}\,\int_{\sqrt 2 \pi \sigma(v - L/2)}^{\sqrt 2 \pi \sigma(v + L/2)}\,{\mathrm e}^{-w^2}\,{\mathrm d}w \,.
$$
For $v \in [- \delta,\, \delta]$, the function $\eta$ can be evaluated as
\begin{eqnarray*}
	\eta(v) &=& \frac{1}{\sqrt \pi}\,\Bigg[ \int_{\mathbb R} {\mathrm e}^{-w^2}\,{\mathrm d}w - \int_{\sqrt 2 \pi \sigma(v - L/2)}^{\sqrt 2 \pi \sigma(v + L/2)}\,{\mathrm e}^{-w^2}\,
	{\mathrm d}w\Bigg]\\ [1ex]
	&=& \frac{1}{\sqrt \pi}\,\Bigg[\int_{-\infty}^{\sqrt 2 \pi \sigma(v - L/2)}\,{\mathrm e}^{-w^2}\,{\mathrm d}w + \int_{\sqrt 2 \pi \sigma(v + L/2)}^{\infty}\,
	{\mathrm e}^{-w^2}\,{\mathrm d}w\Bigg]\\ [1ex]
	&=& \frac{1}{\sqrt \pi}\,\Bigg[\int_{\sqrt 2 \pi \sigma(L/2 - v)}^{\infty}\,{\mathrm e}^{-w^2}\,{\mathrm d}w + \int_{\sqrt 2 \pi \sigma(v + L/2)}^{\infty}\,
	{\mathrm e}^{-w^2}\,{\mathrm d}w\Bigg]\,.
\end{eqnarray*}
By \cite[p.~298, Formula 7.1.13]{abst}, for $x \ge 0$ it holds the inequality
$$
\frac{1}{x + \sqrt{x^2 + 2}}\,{\mathrm e}^{-x^2} \le \int_x^{\infty} {\mathrm e}^{-w^2}\,{\mathrm d}w \le \frac{1}{x + \sqrt{x^2 + 4/\pi}}\,{\mathrm e}^{-x^2}\,,
$$
which can be simplified to
\begin{equation}
	\label{eq:intexp(-t2)}
	\int_x^{\infty} {\mathrm e}^{-w^2}\,{\mathrm d}w \le \frac{1}{2x}\,{\mathrm e}^{-x^2}\,, \quad x > 0\,.
\end{equation}
Therefore, the auxiliary function $\eta(v)$ can be estimated by
$$
\eta(v) < \frac{1}{\sqrt \pi}\,\left( \frac{{\mathrm e}^{-2 \pi^2 \sigma^2 (L/2 - v)^2}}{2 \,\sqrt 2\pi  \sigma (L/2 - v)} + \frac{{\mathrm e}^{-2 \pi^2 \sigma^2 (L/2 + v)^2}}{2 \,\sqrt 2 \pi \sigma (L/2 + v)}\right)\,.
$$
Since the function $\frac{1}{x}\,{\mathrm e}^{- \sigma^2 x^2/2}$ decreases for $x>0$, and  $L/2 - v$, $L/2 +  v \in$ \mbox{$[L/2 - \delta, \, L/2 + \delta]$} by $v \in [- \delta,\, \delta]$ with $0 < \delta <L/2$, we conclude that
$$
\eta(v) < \frac{{\mathrm e}^{-2 \pi^2 \sigma^2 (L/2 - v)^2}}{\sqrt{2 \pi}\,\pi \sigma (L/2 - v)} \le \frac{{\mathrm e}^{-2 \pi^2 \sigma^2 (L/2 - \delta)^2}}{\sqrt{2 \pi}\,\sigma (L \pi - \delta)}\,.
$$
Hence, by \eqref{eq:E1} and \eqref{eq:eta} we receive
\begin{align}
	\label{eq:E1_Gaussian}
	E_{1}(m,\delta,L) &\le \frac{\sqrt{\delta}}{\sqrt{\pi}\,\pi \sigma (L/2 - \delta)}\,{\mathrm e}^{-2 \pi^2\sigma^2 (L/2 - \delta)^2} \,.
\end{align}

Now we examine the truncation error constant \eqref{eq:E2}. Here it holds
\begin{align*}
	\varphi_{\mathrm{Gauss}}^2(\tfrac{m}{L}\big) + L\int_{m/L}^{\infty} \varphi_{\mathrm{Gauss}}^2(t)\,{\mathrm d}t
	=
	\e^{-m^2/(L^2\sigma^2)} + L\sigma \int_{m/(L\sigma)}^{\infty} \e^{-t^2} \,\mathrm dt\,.
\end{align*}
From \eqref{eq:intexp(-t2)} it follows
\begin{align*}
	\e^{-m^2/(L^2\sigma^2)} + L\sigma \int_{m/(L\sigma)}^{\infty} \e^{-t^2} \,\mathrm dt
	\leq
	\frac{2m+L^2\sigma^2}{2m} \,\e^{-m^2/(L^2\sigma^2)} \,.
\end{align*}
Thus, by \eqref{eq:E1} we obtain
\begin{align}
	\label{eq:E2_Gaussian}
	E_2(m,\delta,L) &\le \frac{\sqrt{2L}}{\pi m} \,\sqrt{\frac{2m+L^2\sigma^2}{2m}} \,\e^{-m^2/(2L^2\sigma^2)} \,.
\end{align}
For the special parameter $\sigma = \sqrt{\frac{m}{\pi L\,(L - 2\delta)}}$, both error terms \eqref{eq:E1_Gaussian} and \eqref{eq:E2_Gaussian} have the same exponential decay
such that
\begin{align*}
	E_1(m,\delta,L) &\le \frac{2}{\pi}\,\sqrt{\frac{\delta L}{m\, (L - 2\delta)}}\,{\mathrm e}^{-\pi m (L/2 - \delta)/L} \,, \\
	E_2(m,\delta,L) &\le \frac{1}{\pi}\,\sqrt{\frac{L(2\pi(L-2\delta)+1)/m}{{m\pi (L-2\delta)}}} \,\e^{-\pi m (L/2 - \delta)/L} \,.
\end{align*}
For \mbox{$\delta\in(0,N/2)$} and \mbox{$m\in\N\setminus\{1\}$} it additionally holds
\begin{equation*}
	\sqrt{2\pi(L-2\delta)+1} \leq \sqrt{L}\,\sqrt{2\pi+1} \leq \sqrt{L} (m+1) \,.
\end{equation*}
This completes the proof.
\qedsymbol
\medskip

We remark that Theorem~\ref{Theorem:GaussianWKS} improves the corresponding results in \cite{Q03} and \cite{LZ16} since the actual decay rate could be improved in Theorem~\ref{Theorem:GaussianWKS} from $(m-1)$ to $m$.

\begin{Example}
	\label{ex:visualization_GaussianWKS}
	We aim to visualize the error bound from Theorem \ref{Theorem:GaussianWKS}.
	For a given function \mbox{$f\in {\mathcal B}_{\delta}(\mathbb R)$} with \mbox{$\delta=\tau N \in (0,\,N/2)$} and \mbox{$L= N(1+\lambda)$}, where \mbox{$0<\tau<\frac 12$} and \mbox{$\lambda\geq 0$}, we consider the approximation error
	\begin{equation}
		\label{eq:e_ndelta}
		e_{m,\tau,\lambda}(f) \coloneqq \max_{t\in [-1,\,1]}| f(t) - (R_{\varphi,m} f)(t) | \,.
	\end{equation}
	For $\varphi = \varphi_{\mathrm{Gauss}}$ we show that by \eqref{eq:error_result_Gaussian} it holds $e_{m,\tau,\lambda}(f) \le E_{m,\tau,\lambda}\,\|f\|_{L^2(\mathbb R)}$ where
	\begin{equation}
		\label{eq:E_ndelta_Gauss}
		E_1(m,\delta,L) + E_2(m,\delta,L)
		\le E_{m,\tau,\lambda}
		\coloneqq \frac{2\sqrt{\pi \delta L} + L(m+1)/\sqrt m}{\pi\,\sqrt{m\pi(L- 2 \delta)}}\,{\mathrm e}^{-\pi m (L/2 - \delta)/L} \,
	\end{equation}
	with \mbox{$\sigma = \sqrt{\frac{m}{\pi L\,(L - 2\delta)}}$}.
	The error \eqref{eq:e_ndelta} shall here be approximated by evaluating a given function $f$ and its approximation $R_{\varphi,m} f$ at \mbox{$S=10^5$} equidistant points \mbox{ $t_s\in [-1,\,1]$}, \mbox{$s=1,\dots,S$}.
	By the definition of the regularized Shannon sampling formula in \eqref{eq:Rmf} it can be seen that for $t \in [-1, 1]$ we have
	$$(R_{\varphi,m} f)(t) = \sum_{\ell = -L-m}^{L+m} f(\tfrac \ell L) \,\psi(t - \tfrac \ell L) \,.$$
	Here we study the function \mbox{$f(t) = \sqrt{2 \delta} \,\mathrm{sinc}(2 \delta \pi t)$}, \mbox{$t \in \mathbb R$}, such that it holds \mbox{$\|f\|_{L^2(\mathbb R)}=1$}.
	We fix \mbox{$N=128$} and consider the evolution for different values \mbox{$m \in \mathbb N \setminus \{1\}$}, i.\,e., we are still free to make a choice for the parameters $\tau$ and $\lambda$.
	In a first experiment we fix \mbox{$\lambda=1$} and choose different values for $\tau<\frac 12$, namely we consider \mbox{$\tau \in \{\nicefrac{1}{20},\,\nicefrac{1}{10},\,\nicefrac{1}{4},\,\nicefrac{1}{3},\,\nicefrac{9}{20}\}$}.
	The corresponding results are depicted in Figure~\ref{fig:err_const_Gaussian}~(a).
	We recognize that the smaller the factor~$\tau$ can be chosen, the better the error results are.
	As a second experiment we fix \mbox{$\tau=\frac 13$}, but now choose different \mbox{$\lambda\in\{0,0.5,1,2\}$}.
	The associated results are displayed in Figure~\ref{fig:err_const_Gaussian}~(b).
	It can clearly be seen that the higher the oversampling parameter $\lambda$ is chosen, the better the error results get.
	We remark that for larger choices of $N$, the line plots in Figure~\ref{fig:err_const_Gaussian} would only be shifted slightly upwards, such that for all $N$ we receive almost the same error results.
	\begin{figure}[ht]
		\centering
		\captionsetup[subfigure]{justification=centering}
		\begin{subfigure}[t]{0.49\textwidth}
			\includegraphics[width=\textwidth]{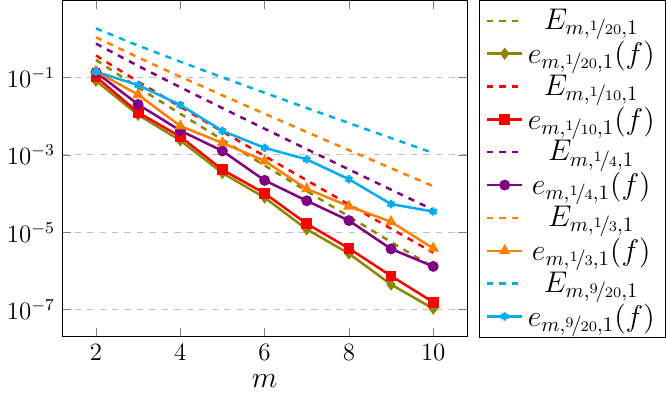}
			\caption{$\lambda=1$ and various $\tau<\frac 12$}
		\end{subfigure}
		\begin{subfigure}[t]{0.49\textwidth}
			\includegraphics[width=\textwidth]{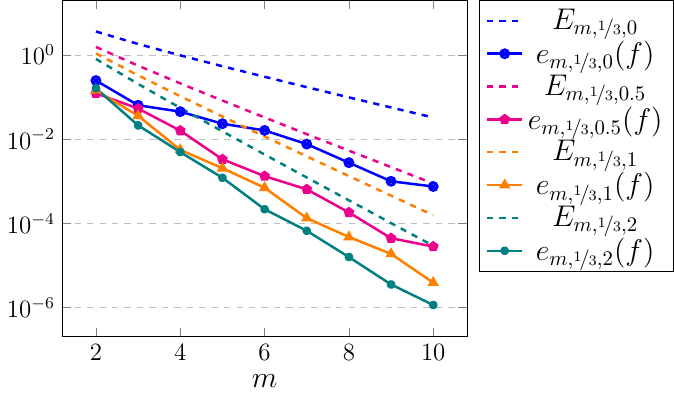}
			\caption{$\tau=\frac 13$ and various $\lambda\geq 0$}
		\end{subfigure}
		\caption{Maximum approximation error \eqref{eq:e_ndelta} and error constant \eqref{eq:E_ndelta_Gauss} using \mbox{$\varphi_{\mathrm{Gauss}}$} in \eqref{eq:varphiGauss} and \mbox{$\sigma = \sqrt{\frac{m}{\pi L\,(L - 2\delta)}}$} for the function \mbox{$f(x) = \sqrt{2\delta} \,\mathrm{sinc}(2 \delta \pi x)$} with \mbox{$N=128$}, \mbox{$m\in\{2, 3, \ldots, 10\}$}, as well as \mbox{$\tau \in \{\nicefrac{1}{20},\,\nicefrac{1}{10},\,\nicefrac{1}{4},\,\nicefrac{1}{3},\,\nicefrac{9}{20}\}$}, \mbox{$\delta = \tau N$}, and \mbox{$\lambda\in\{0,0.5,1,2\}$}, respectively.}
		\label{fig:err_const_Gaussian}
	\end{figure}
\end{Example}

Now we show that for the regularized Shannon sampling formula with the Gaussian window function~\eqref{eq:varphiGauss} the uniform perturbation error \eqref{eq:perturb_error} only grows as $\O{\sqrt{m}}$.
We remark that a similar result can also be found in \cite{Q05}.

\begin{Theorem}
	\label{Theorem:robustness_Gaussian}
	Let $f\in {\mathcal B}_{\delta}(\mathbb R)$ with $\delta = \tau N$, $\tau\in (0,\,1/2)$, $N \in \mathbb N$, \mbox{$L= N(1+\lambda)$} with \mbox{$\lambda\geq 0$} \\and \mbox{$m \in {\mathbb N}\setminus \{1\}$} be given.
	Further let \mbox{$R_{\mathrm{Gauss},m}{\tilde f}$} be as in \eqref{eq:Rmf_erroneous} with the noisy samples
	${\tilde f}_{\ell} = f\big(\frac{\ell}{L}\big) + \varepsilon_{\ell}$,
	where $|\varepsilon_{\ell}| \le \varepsilon$ for all $\ell\in\Z$ and $\varepsilon>0$.\\
	Then the regularized Shannon sampling formula \eqref{eq:Rmf} with the Gaussian window function~\eqref{eq:varphiGauss} and $\sigma = \sqrt{\frac{m}{\pi L\,(L - 2\delta)}}$ satisfies
	\begin{align}
		\label{eq:result_robustness_Gauss}
		\| R_{\mathrm{Gauss},m}{\tilde f} - R_{\mathrm{Gauss},m}f \|_{C_0(\mathbb R)}
		\leq \varepsilon \left(2 + \sqrt{\frac{2 + 2\lambda}{\lambda + 1 - 2 \tau}}\, \sqrt m\right) \,.
	\end{align} 	
\end{Theorem}

{\emph Proof.} By Theorem \ref{Theorem:robustness} we only have to compute $\hat{\varphi}_{\mathrm{Gauss}}(0)$ for the Gaussian window function \eqref{eq:varphiGauss}.
By \eqref{eq:Gaussian_FT} we recognize that
\begin{align*}
	\hat{\varphi}_{\mathrm{Gauss}}(0) &= \sqrt{2\pi}\,\sigma = \sqrt{\frac{2 m}{L\,(L - 2\delta)}} = \frac{1}{L}\,\sqrt{\frac{2 + 2\lambda}{\lambda + 1 - 2 \tau}}\, \sqrt m
\end{align*}
such that \eqref{eq:result_robustness} yields the assertion.
\qedsymbol
%\medskip

\begin{Example}
\label{ex:perturbation_Gaussian}
Now we visualize the error bound from Theorem~\ref{Theorem:robustness_Gaussian}.
Similar to Example~\ref{ex:visualization_GaussianWKS}, we consider the perturbation error
\begin{equation}
	\label{eq:etilde_ndelta}
	\tilde e_{m,\tau,\lambda}(f) \coloneqq \max_{t\in [-1,\,1]}| (R_{\varphi,m} \tilde f)(t) - (R_{\varphi,m} f)(t) | \,.
\end{equation}
For $\varphi = \varphi_{\mathrm{Gauss}}$ we show that by \eqref{eq:result_robustness_Gauss} it holds
$\tilde e_{m,\tau,\lambda}(f) \le \tilde E_{m,\tau,\lambda}$, where
\begin{equation}
	\label{eq:Etilde_ndelta_Gauss}
	\tilde E_{m,\tau,\lambda}
	\coloneqq \varepsilon \left(2 + \sqrt{\frac{2 + 2\lambda}{\lambda + 1 - 2 \tau}}\, \sqrt m\right) \,.
\end{equation}
We conduct the same experiments as in Example~\ref{ex:visualization_GaussianWKS} and introduce a maximum perturbation of \mbox{$\varepsilon=10^{-3}$} as well as uniformly distributed random numbers $\varepsilon_{\ell}$ in $(-\varepsilon,\varepsilon)$.
Due to the randomness we perform the experiments one hundred times and then take the maximum error over all runs.
The associated results are displayed in Figure~\ref{fig:err_perturbation_Gaussian}.
\begin{figure}[ht]
	\centering
	\captionsetup[subfigure]{justification=centering}
	\begin{subfigure}[t]{0.49\textwidth}
		\includegraphics[width=\textwidth]{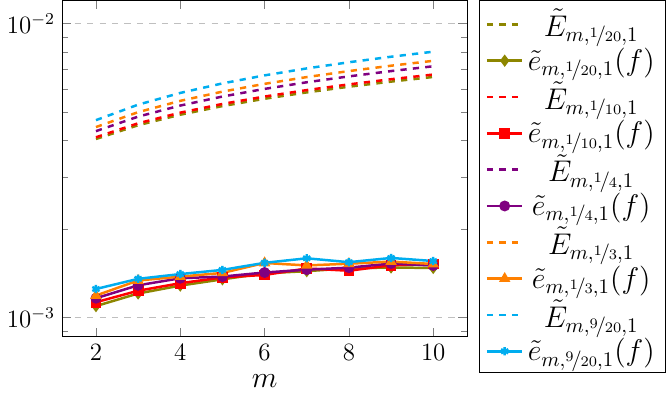}
		\caption{$\lambda=1$ and various $\tau<\frac 12$}
	\end{subfigure}
	\begin{subfigure}[t]{0.49\textwidth}
		\includegraphics[width=\textwidth]{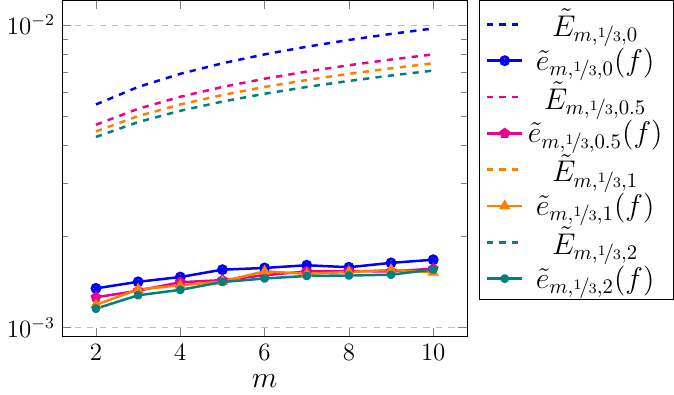}
		\caption{$\tau=\frac 13$ and various $\lambda\geq 0$}
	\end{subfigure}
	\caption{Maximum perturbation error \eqref{eq:etilde_ndelta} over 100 runs and error constant \eqref{eq:Etilde_ndelta_Gauss} using \mbox{$\varphi_{\mathrm{Gauss}}$} in \eqref{eq:varphiGauss} and \mbox{$\sigma = \sqrt{\frac{m}{\pi L\,(L - 2\delta)}}$} for the function \mbox{$f(x) = \sqrt{2\delta} \,\mathrm{sinc}(2 \delta \pi x)$} with \mbox{$\varepsilon=10^{-3}$}, \mbox{$N=128$}, \mbox{$m\in\{2, 3, \ldots, 10\}$}, as well as \mbox{$\tau \in \{\nicefrac{1}{20},\,\nicefrac{1}{10},\,\nicefrac{1}{4},\,\nicefrac{1}{3},\,\nicefrac{9}{20}\}$}, \mbox{$\delta = \tau N$}, and \mbox{$\lambda\in\{0,0.5,1,2\}$}, respectively.}
	\label{fig:err_perturbation_Gaussian}
\end{figure}
\end{Example}

\section{B--spline regularized Shannon sampling formula\label{sec:Bspline}}

Now we consider the modified $\mathrm B$--spline window function~\eqref{eq:varphiB}
with $s,\,m \in \mathbb N \setminus \{1\}$ and $L= N\,(1 + \lambda)$, $\lambda \ge 0$, where $M_{2s}$ denotes the
centered cardinal $\mathrm B$--spline of even order~$2s$. Note that \eqref{eq:varphiB} is supported on $\big[-\frac{m}{L},\, \frac{m}{L}\big]$.

\begin{Lemma}
\label{Lemma:M2m(0)}
For the value $M_{2s}(0)$, $s\in \mathbb N$, it holds the formula
\begin{equation}
\label{eq:M2m(0)}
M_{2s}(0) = \frac{1}{(2s-1)!}\,\sum_{j=0}^{s-1} (-1)^j\,{2s \choose j}\,(s-j)^{2s-1}\,.
\end{equation}
The sequence $\big(\sqrt{2s}\,M_{2s}(0)\big)_{s=1}^{\infty}$ has the limit
\begin{equation}
\label{eq:limM2m(0)}
\lim_{s \to \infty} \sqrt{2s}\,M_{2s}(0) = \sqrt{\frac{6}{\pi}} \approx 1.3820\,.
\end{equation}
\end{Lemma}

\emph{Proof.} By inverse Fourier transform of ${\hat \varphi}_{\mathrm B}$ it holds
$$
\varphi_{\mathrm B}(x) = \int_{\mathbb R} {\hat \varphi}_{\mathrm B}(v)\,{\mathrm e}^{2 \pi {\mathrm i}v x}\, {\mathrm d}v\,, \quad x \in \mathbb R\,.
$$
Hence, for $x = 0$ it follows that
$$
M_{2s}(0) = \int_{\mathbb R} \big({\mathrm{sinc}} (\pi v)\big)^{2s}\,{\mathrm d}v = \frac{2}{\pi}\,\int_0^{\infty} \big({\mathrm{sinc}} \,w\big)^{2s}\,{\mathrm d}w\,.
$$
The above integral can be determined in explicit form (see \cite{MeRo65}, \cite[p.~20, 5.12]{Ob90} or \cite[(4.1.12)]{chui92}) as
$$
\int_0^{\infty} \big({\mathrm{sinc}} \,w\big)^{2s}\,{\mathrm d}w = \frac{\pi}{2\,(2s-1)!}\,\sum_{j=0}^{s-1} (-1)^j\,{2s \choose j}\,(s-j)^{2s-1}
$$
such that \eqref{eq:M2m(0)} is shown. Especially, it holds $M_2(0) = 1$, $M_4(0) = \frac{2}{3}$, $M_6(0) = \frac{11}{20}$, $M_8(0) = \frac{151}{315}$,
$M_{10}(0) = \frac{15619}{36288}$, and $M_{12}(0) = \frac{655177}{1663200}$.
A table with the decimal values of $M_{2s}(0)$ for $m =15,\,\ldots,\,50,$ can be found in \cite{MeRo65}.
For example, it holds $M_{100}(0) \approx 0.137990$.

By \cite[(3.6)]{UAE}, there exists the pointwise limit
$$
\lim_{s \to \infty} \sqrt{\frac{s}{6}}\, M_{2s}\left(\sqrt{\frac{s}{6}}\,x \right) = \frac{1}{\sqrt{2\pi}}\,{\mathrm e}^{-x^2/2}
$$
such that for $x = 0$ we obtain \eqref{eq:limM2m(0)}.
\qedsymbol
\medskip

\begin{Remark}
By numerical computations we can see that the sequence  $\big(\sqrt{2s}\,M_{2s}(0)\big)_{s=2}^{50}$ increases monotonously, see Figure~\ref{fig:assertion_bspline}.
For large $s$ we can use the asymptotic expansion
$$
\sqrt{2s}\,M_{2s}(0) \approx \sqrt{\frac{6}{\pi}}\,\left[1 - \frac{3}{40\,s} - \frac{13}{4480\,s^2} + \frac{27}{25600\,s^3} + \frac{52791}{63078400\,s^4} + \frac{482427}{2129920000\,s^5}\right]
$$
(see~\cite{MeRo65}) such that the whole sequence $\big(\sqrt{2s}\,M_{2s}(0)\big)_{s=2}^{\infty}$ increases monotonously.
Hence, for $s \in {\mathbb N}\setminus \{1\}$ it holds
\begin{equation}
\label{eq:estM2m(0)}
\frac{4}{3} \le \sqrt{2s}\,M_{2s}(0) < \sqrt{\frac{6}{\pi}} \,, \quad s \in \mathbb N \setminus \{1\} \,.
\end{equation}
\begin{figure}[ht]
	\centering
	\captionsetup[subfigure]{justification=centering}
	\begin{subfigure}[t]{0.4\textwidth}
		\includegraphics[width=\textwidth]{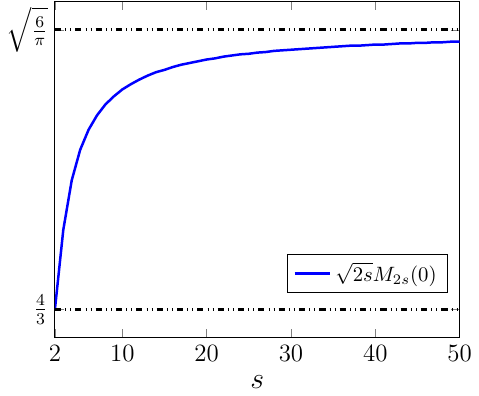}
	\end{subfigure}
	\caption{The sequence $\big(\sqrt{2s}\,M_{2s}(0)\big)_{s=2}^{50}$.}
	\label{fig:assertion_bspline}
\end{figure}
\end{Remark}

Now we show that for the regularized Shannon sampling formula~\eqref{eq:Rmf} with $\mathrm B$--spline window function \eqref{eq:varphiB} the uniform approximation error \eqref{eq:Error_RWKS} decays exponentially with respect to $m$.

\begin{Theorem}
\label{Theorem:bsplineWKS}
Let $f\in {\mathcal B}_{\delta}(\mathbb R)$ with $\delta = \tau N$, $\tau \in (0,\,1/2)$, $N\in \mathbb N$, $L = N\,(1 + \lambda)$, $\lambda \geq 0$, and $m \in \mathbb N \setminus \{1\}$ be given. Assume that
\begin{equation}
	\label{eq:condtau}
	\frac{\tau}{1 + \lambda} < \frac{1}{2} - \frac{1}{\pi}\,.
\end{equation}
Then the regularized Shannon sampling formula \eqref{eq:Rmf} with the $\mathrm B$--spline window function~\eqref{eq:varphiB} and $s = \left\lceil \frac{m+1}{2} \right\rceil$ satisfies the error estimate
\begin{align}
\label{eq:error_result_bspline}
\| f - R_{{\mathrm B},m} f \|_{C_0(\mathbb R)} \le \frac{3\,\sqrt{\delta s}}{(2s-1)\,\pi}\, \e^{-m\,\left( \ln\left(\pi m\,(1 + \lambda - 2 \tau)\right) - \ln \left(2s(1 + \lambda)\right) \right)} \, \| f \|_{L^2(\mathbb R)}\,.
\end{align}
\end{Theorem}

{\emph Proof.} By Theorem \ref{Theorem:errorRmf} we only have to estimate the regularization error constant \eqref{eq:E1},
since it holds
\mbox{$\varphi_{\mathrm B}(x) = \varphi_{\mathrm B}(x)\,{\mathbf 1}_{[-m/L,\,m/L]}(x)$} for all $x \in \mathbb R$ and therefore the truncation error constant \eqref{eq:E2} vanishes for the $\mathrm B$--spline window function \eqref{eq:varphiB} by Remark \ref{Remark:simpleerror}.

By \cite{PPST18, PT21a} it holds
\begin{align}
\label{eq:bspline_FT}
	{\hat \varphi_{\mathrm{B}}}(v) &= \frac{m}{sL\, M_{2s}(0)}\,\left(\mathrm{sinc} \,\frac{\pi vm}{sL}\right)^{2s}\,, \quad v \in \mathbb R \,,
\end{align}
such that
the auxiliary function \eqref{eq:eta} is given by
$$
\eta(v) = {\mathbf 1}_{[-\delta,\, \delta]}(v) - \frac{m}{sL\,M_{2s}(0)}\,\int_{v-L/2}^{v+L/2} \left(\mathrm{sinc} \,\frac{\pi um}{sL}\right)^{2s}\,{\mathrm d}u \,, \quad v \in \mathbb R \,.
$$
By inverse Fourier transform \eqref{eq:inverse_Fouriertrafo} we have
\begin{equation}
\label{eq:estintsinc}
1 = \varphi_{\mathrm B}(0) = \int_{\mathbb R} {\hat \varphi}_{\mathrm B}(v) \,{\mathrm d}v = \frac{m}{sL\,M_{2s}(0)}\, \int_{\mathbb R} \left(\mathrm{sinc}\,
\frac{\pi vm}{sL}\right)^{2s} \,{\mathrm d}v \,.
\end{equation}
Then, for $v \in [- \delta, \, \delta]$, the function $\eta$ can be determined by \eqref{eq:estintsinc} in the following form
\begin{eqnarray*}
	\eta(v)
	&=& \frac{m}{sL\,M_{2s}(0)}\,\left[\int_{\mathbb R} \left(\mathrm{sinc} \,\frac{\pi um}{sL}\right)^{2s}\,{\mathrm d}u - \int_{v- L/2}^{v+L/2} \left(\mathrm{sinc} \,\frac{\pi um}{sL}\right)^{2s}\,{\mathrm d}u\right]\\
	&=& \frac{m}{sL\,M_{2s}(0)}\,\left[\int_{L/2-v}^{\infty} \left(\mathrm{sinc} \,\frac{\pi um}{sL}\right)^{2s}\,{\mathrm d}u + \int_{v + L/2}^{\infty} \big(\mathrm{sinc}
	\,\frac{\pi um}{sL}\big)^{2s}\,{\mathrm d}u\right]\,.
\end{eqnarray*}
Applying the simple estimates
\begin{eqnarray*}
	\int_{L/2-v}^{\infty} \left(\mathrm{sinc} \,\frac{\pi um}{sL}\right)^{2s}\,{\mathrm d}u &\le& \frac{s^{2s}\,L^{2s}}{m^{2s}\,\pi^{2s}}\, \int_{L/2-v}^{\infty} u^{-2s}\,{\mathrm d}u = \frac{s^{2s}\,L^{2s}}{(2s-1)\,
		m^{2s}\,\pi^{2s}\,(L/2-v)^{2s-1}}\,,\\
	\int_{v + L/2}^{\infty} \left(\mathrm{sinc} \,\frac{\pi um}{sL}\right)^{2s}\,{\mathrm d}u &\le& \frac{s^{2s}\,L^{2s}}{m^{2s}\,\pi^{2s}}\, \int_{v + L/2}^{\infty} u^{-2s}\,{\mathrm d}u = \frac{s^{2s}\,L^{2s}}{(2s-1)\,
		m^{2s}\,\pi^{2s}\,(v + L/2)^{2s-1}}\,,
\end{eqnarray*}
the function $\eta$ can be estimated for $v \in [- \delta, \delta]$ by
$$
\eta(v) \le \frac{s^{2s-1}\,L^{2s-1}}{(2s-1)\,m^{2s-1}\,\pi^{2s}\,M_{2s}(0)}\,\left[\frac{1}{(L/2 - v)^{2s-1}} + \frac{1}{(L/2 + v)^{2s-1}}\right]\,.
$$
By $v \in [- \delta,\,\delta]$ with $0 < \delta < N/2 \le L/2$, it holds $L/2 - v$, $L/2 + v \in [L/2 - \delta,\,L/2 + \delta]$. Since the function $x^{1-2s}$ decreases for
$x > 0$, we conclude that
$$
\max_{v \in [-\delta,\delta]} |\eta(v)| \leq \frac{2\,s^{2s-1}\,L^{2s-1}}{(2s-1)\,m^{2s-1}\,\pi^{2s}\,M_{2s}(0)\,(L/2 - \delta)^{2s-1}} \,.
$$
Hence, by \eqref{eq:E1}, \eqref{eq:eta} and \eqref{eq:estM2m(0)} we receive
\begin{align}
\label{eq:E1_B}
	E_1(m,\delta,L)
	&\le
	\frac{2\,\sqrt{2 \delta}}{(2s-1)\,\pi\,M_{2s}(0)} \left(\frac{2 sL}{\pi mL - 2 \pi m\delta}\right)^{2s-1}
	\le
	\frac{3\,\sqrt{\delta s}}{(2s-1)\,\pi}\, \left(\frac{2 sL}{\pi mL - 2 \pi m\delta}\right)^{2s-1} \,.
\end{align}
To guarantee convergence of this result, we have to satisfy
\begin{equation}
\label{eq:cond}
	\frac{2 sL}{\pi mL - 2 \pi m\delta} = \frac{2s(1 + \lambda)}{\pi m\,(1 + \lambda - 2 \tau)} \eqqcolon c < 1 \,.
\end{equation}
By means of logarithmic laws we recognize that
$c^{\,2s-1} = \e^{\ln (c^{2s-1})} = \e^{(2s-1)\,\ln c}$.
Thus, the condition $c < 1$ yields $\ln c<0$ and therefore an exponential decay of \eqref{eq:E1_B} with respect to $(2s-1)$.
Thereby, we need on the one hand that $\ln c$ is as small as possible, which is equivalent to choosing $s$ as small as possible.
On the other hand, we aim achieving a decay rate of at least~$m$, i.\,e., we need to fulfill $2s-1\geq m$.
These two conditions
can now be used to pick the optimal parameter $s\in \mathbb N$ in the form
$s = \left\lceil \frac{m+1}{2} \right\rceil$.
Then \eqref{eq:cond} holds, if \eqref{eq:condtau} is fulfilled, and
\begin{align*}
	c^{\,2s-1}
	&=
	\e^{(2s-1)\,\ln c}
	=
	\e^{(2s-1)\,\left(\ln \left(2s(1 + \lambda)\right) - \ln\left(\pi m\,(1 + \lambda - 2 \tau)\right)\right)} \\
	&=
	\e^{-(2s-1)\,\left( \ln\left(\pi m\,(1 + \lambda - 2 \tau)\right) - \ln \left(2s(1 + \lambda)\right) \right)}
	\leq
	\e^{-m\,\left( \ln\left(\pi m\,(1 + \lambda - 2 \tau)\right) - \ln \left(2s(1 + \lambda)\right) \right)}
\end{align*}
yields the assertion.
We remark that it holds $\pi m\,(1 + \lambda - 2 \tau) > 2s(1 + \lambda)$ since $c<1$.
\qedsymbol

\begin{Example}
\label{ex:visualization_bsplineWKS}
Analogous to Example~\ref{ex:visualization_GaussianWKS}, we now visualize the error bound from Theorem \ref{Theorem:bsplineWKS}, i.\,e.,
for $\varphi = \varphi_{\mathrm{B}}$ we show that for the approximation error \eqref{eq:e_ndelta} it holds by \eqref{eq:error_result_bspline} that $e_{m,\tau,\lambda}(f) \le E_{m,\tau,\lambda}\,\|f\|_{L^2(\mathbb R)}$, where
\begin{equation}
	\label{eq:E_ndelta_B}
	E_1(m,\delta,L) \le E_{m,\tau,\lambda} \coloneqq \frac{3\,\sqrt{\delta s}}{(2s-1)\,\pi}\, \e^{-m\,\left( \ln\left(\pi m\,(1 + \lambda - 2 \tau)\right) - \ln \left(2s(1 + \lambda)\right) \right)}
\end{equation}
with $s = \left\lceil \frac{m+1}{2} \right\rceil$.
Additionally, we now have to observe the condition \eqref{eq:condtau}.
For the first experiment in Example~\ref{ex:visualization_GaussianWKS} with $\lambda=1$ this leads to $\tau<1-\frac{2}{\pi}\approx 0.3634$,
while in the second experiment we fixed $\tau = \frac 13$ and therefore have to satisfy $\lambda>\frac{2\pi}{3\pi-6}-1\approx 0.8346$.
Thus, only in these settings the requirements of Theorem~\ref{Theorem:bsplineWKS} are fulfilled, and therefore only those error bounds are plotted in Figure~\ref{fig:err_const_bspline} while the approximation error \eqref{eq:e_ndelta} is computed for all constellations of parameters as given in Example~\ref{ex:visualization_GaussianWKS}.
We recognize that we have almost the same behavior as in Figure~\ref{fig:err_const_Gaussian}, which means that there is hardly any improvement using the $\mathrm{B}$--spline window function in comparison to the well-studied Gaussian window function.
\begin{figure}[ht]
	\centering
	\captionsetup[subfigure]{justification=centering}
	\begin{subfigure}[t]{0.49\textwidth}
		\includegraphics[width=\textwidth]{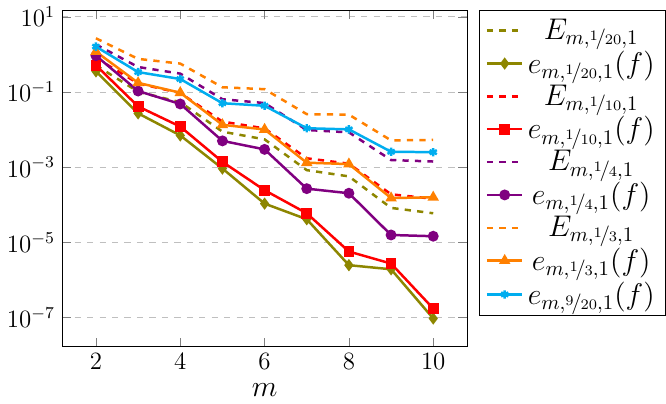}
		\caption{$\lambda=1$ and various $\tau<\frac 12$}
	\end{subfigure}
	\begin{subfigure}[t]{0.49\textwidth}
		\includegraphics[width=\textwidth]{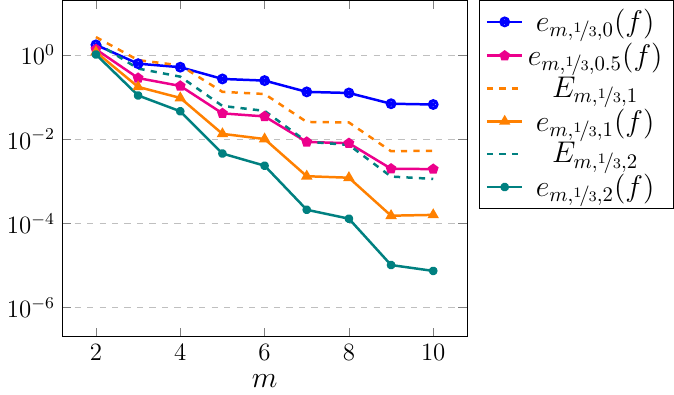}
		\caption{$\tau=\frac 13$ and various $\lambda\geq 0$}
	\end{subfigure}
	\caption{Maximum approximation error \eqref{eq:e_ndelta} and error constant \eqref{eq:E_ndelta_B} using \mbox{$\varphi_{\mathrm{B}}$} in \eqref{eq:varphiB} and \mbox{$s = \left\lceil \frac{m+1}{2} \right\rceil$} for the function \mbox{$f(x) = \sqrt{2\delta} \,\mathrm{sinc}(2 \delta \pi x)$} with \mbox{$N=128$}, \mbox{$m\in\{2, 3, \ldots, 10\}$}, as well as \mbox{$\tau \in \{\nicefrac{1}{20},\,\nicefrac{1}{10},\,\nicefrac{1}{4},\,\nicefrac{1}{3},\,\nicefrac{9}{20}\}$}, \mbox{$\delta = \tau N$}, and \mbox{$\lambda\in\{0,0.5,1,2\}$}, respectively.
		\label{fig:err_const_bspline}}
\end{figure}
\end{Example}

Now we show that for the regularized Shannon sampling formula  with the $\mathrm B$--spline window function \eqref{eq:varphiB} the uniform perturbation error \eqref{eq:result_robustness}
only grows as $\O{\sqrt{m}}$.

\begin{Theorem}
\label{Theorem:robustness_B}
Let $f\in {\mathcal B}_{\delta}(\mathbb R)$ with $\delta = \tau N$, $\tau\in (0,\,1/2)$, $N \in \mathbb N$, \mbox{$L= N(1+\lambda)$} with \mbox{$\lambda\geq 0$} and \mbox{$m \in {\mathbb N}\setminus \{1\}$} be given.
Let $s \in \mathbb N$ be defined by $s = \left\lceil \frac{m+1}{2} \right\rceil$.
Further let \mbox{$R_{\mathrm{B},m}{\tilde f}$} be as in \eqref{eq:Rmf_erroneous} with the noisy samples
${\tilde f}_{\ell} = f\big(\frac{\ell}{L}\big) + \varepsilon_{\ell}$,
where $|\varepsilon_{\ell}| \le \varepsilon$ for all $\ell\in\Z$ and $\varepsilon>0$.\\
Then the regularized Shannon sampling formula \eqref{eq:Rmf} with the $\mathrm B$--spline window function~\eqref{eq:varphiB} and $s = \left\lceil \frac{m+1}{2} \right\rceil$ satisfies
\begin{align}
\label{eq:result_robustness_B}
	\| R_{\mathrm{B},m}{\tilde f} - R_{\mathrm{B},m}f \|_{C_0(\mathbb R)}
	\leq \varepsilon \left( 2+\frac{3}{2} \,\sqrt{m}\right) \,.
\end{align} 	
\end{Theorem}

{\emph Proof.} By Theorem \ref{Theorem:robustness} we only have to compute $\hat{\varphi}_{\mathrm B}(0)$ for the $\mathrm B$--spline window function \eqref{eq:varphiB}.
By \eqref{eq:bspline_FT} and \eqref{eq:estM2m(0)} we recognize that
\begin{align*}
\hat{\varphi}_{\mathrm B}(0) &= \frac{m}{sL\, M_{2s}(0)} \le \frac{m}{sL} \frac{3 \,\sqrt 2}{4}\,\sqrt s \,.
\end{align*}
Due to $s = \left\lceil \frac{m+1}{2} \right\rceil$ it holds $\sqrt s \geq \frac{\sqrt{m}}{\sqrt{2}}$
such that \eqref{eq:result_robustness} yields the assertion \eqref{eq:result_robustness_B}.
\qedsymbol
\medskip

Similar to Example~\ref{ex:perturbation_Gaussian}, one can visualize the error bound from Theorem~\ref{Theorem:robustness_B}, which leads to results analogous to Figure~\ref{fig:err_perturbation_Gaussian}.

\section{sinh-type regularized Shannon sampling formula\label{sec:sinh}}

We consider the $\sinh$-type window function~\eqref{eq:varphisinh}
and show that in this case the uniform approximation error~\eqref{eq:Error_RWKS} for the regularized Shannon sampling formula~\eqref{eq:Rmf} decays exponentially with respect to $m$.

\begin{Theorem}
\label{Thm:errorpsisinh}
Let $f\in {\mathcal B}_{\delta}(\mathbb R)$ with $\delta = \tau N$, $\tau \in (0,\,1/2)$, $N\in \mathbb N$, $L = N\,(1 + \lambda)$, $\lambda \geq 0$, and
$m \in \mathbb N \setminus \{1\}$ be given.\\
Then the regularized Shannon sampling formula \eqref{eq:Rmf} with the $\sinh$-type window function~\eqref{eq:varphisinh} and $\beta = \frac{\pi m\,(1 + \lambda + 2 \tau)}{1 + \lambda}$
satisfies the error estimate
\begin{align*}
	\| f - R_{\sinh,m} f \|_{C_0(\mathbb R)}
	\le
	\bigg(\frac{\sqrt{\beta\,\pi\delta}}{(1 - 2\,{\mathrm e}^{-\beta})(1-w_0^2)^{1/4}}\,
	\,
	{\mathrm e}^{-\beta\big(1-\sqrt{1-w_0^2}\big)} + \frac{2\,\sqrt{2 \delta}}{1 - {\mathrm e}^{-2\beta}}\,{\mathrm e}^{-\beta}\bigg)\,
	\| f \|_{L^2(\mathbb R)} \,,
\end{align*}
where $w_0 = \frac{1 + \lambda - 2\tau}{1 + \lambda + 2\tau} \in (0,\,1)$.
Further, the regularized Shannon sampling formula \eqref{eq:Rmf} with the $\sinh$-type window function \eqref{eq:varphisinh} and $\beta = \frac{\pi m\,(1 + \lambda - 2 \tau)}{1 + \lambda}$ fulfills
\begin{align*}
	\| f - R_{\sinh,m} f \|_{C_0(\mathbb R)}
	\le 3\,\sqrt{2\delta} \, \e^{-\beta} \, \| f \|_{L^2(\mathbb R)} \,.
\end{align*}
\end{Theorem}

\emph{Proof.} By Theorem \ref{Theorem:errorRmf} we only have to estimate the regularization error constant \eqref{eq:E1},
since it holds
\mbox{$\varphi_{\sinh}(x) = \varphi_{\sinh}(x)\,{\mathbf 1}_{[-m/L,\,m/L]}(x)$}
for all $x \in \mathbb R$ and therefore the truncation error constant \eqref{eq:E2} vanishes for the $\sinh$-type window function \eqref{eq:varphisinh} by Remark \ref{Remark:simpleerror}.

By \cite[p.~38, 7.58]{Ob90} or \cite{PT21} it holds
\begin{align}
\label{eq:sinh_FT}
	{\hat \varphi_{\sinh}}(v)
	&= \frac{\pi m \beta}{L\, \sinh \beta} \cdot \left\{\begin{array}{ll} (w^2 - \beta^2)^{-1/2}\,J_1\big(\sqrt{w^2 - \beta^2}\big) & \quad w\in \mathbb R \setminus \{- \beta, \, \beta\}\,,\\
		1/2 & \quad w = \pm \beta,
	\end{array} \right.
\end{align}
where $w \coloneqq 2\pi m v/L$ denotes a scaled frequency.
Thus, the auxiliary function \eqref{eq:eta} is given by
$$
\eta(v) = {\mathbf 1}_{[-\delta,\, \delta]}(v) - \frac{m\beta}{2 L\,\sinh \beta}\,\int_{v-L/2}^{v+L/2} \frac{J_1\big(2 \pi \,\sqrt{m^2u^2/L^2 - s^2(1+ 2 \lambda)^2/(2 + 2\lambda)^2}\big)}
{\sqrt{m^2u^2/L^2 - s^2(1+ 2 \lambda)^2/(2 + 2\lambda)^2}}\,{\mathrm d}u \,, \ v \in \mathbb R \,.
$$
Substituting $u = \frac{sL\,(1+ 2\lambda)}{m(2 + 2\lambda)}\,w$, we obtain for $v \in [- \delta, \, \delta]$ that
\begin{align}
\label{eq:eta(v)}
	\eta(v) = 1 - \frac{\beta}{2\,\sinh \beta}\,\int_{-w_1(-v)}^{w_1(v)} \frac{J_1\big(\beta\, \sqrt{w^2-1}\big)}{\sqrt{w^2-1}}\,{\mathrm d}w
\end{align}
with
\begin{align}
\label{eq:w1(v)}
	w_1(v) \coloneqq \frac{m(v + L/2)(2+2\lambda)}{sL\,(1+2\lambda)} >0\,, \quad v \in [- \delta, \, \delta] \,.
\end{align}
Since the integrand of \eqref{eq:eta(v)} behaves differently for \mbox{$w \in [-1, 1]$} and \mbox{$w \in \mathbb R \setminus (-1, 1)$} we have to distinguish between the cases \mbox{$w_1(v)\leq 1$} and \mbox{$w_1(v)\geq 1$} for all \mbox{$v \in [-\delta,\, \delta]$}.
By definition $w_1(v)$ is linear and monotonously increasing.
Thus, we have
\mbox{$\min \{w_1(v):\,v\in [-\delta,\,\delta]\} = w_1(-\delta)$} and
\mbox{$\max \{w_1(v):\,v\in [-\delta,\,\delta]\} = w_1(\delta)$}.
This fact can now be used to choose an optimal parameter \mbox{$s = s(m,\tau,\lambda) > 0$} such that either \mbox{$w_1(\delta)\leq 1$} or \mbox{$w_1(-\delta)\geq 1$} is fulfilled.

\begin{figure}[t]
	\centering
	\captionsetup[subfigure]{justification=centering}
	\begin{subfigure}{0.4\textwidth}
		\includegraphics[width=\textwidth]{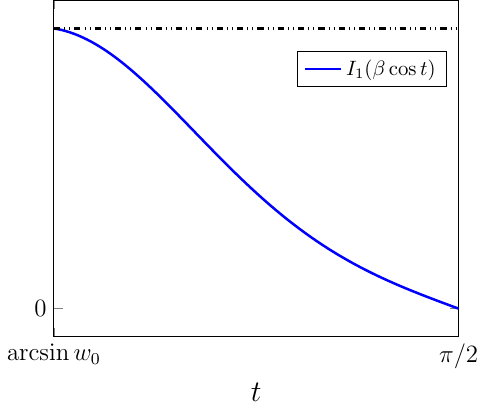}
	\end{subfigure}
	\caption{The integrand $I_1(\beta\,\cos t)$ on the interval $[{\arcsin\,w_0},\,{\pi/2}]$.}
	\label{fig:integrand_i1}
\end{figure}

\paragraph{Case 1 ($\boldsymbol{w_1(\delta)\leq 1}$):}
Note that by \cite[6.681--3]{GR80} and \cite[10.2.13]{abst} as well as $J_1({\mathrm i}\,z) = {\mathrm i}\,I_1(z)$ for $z \in \mathbb C$ it holds
\begin{align}
\label{eq:integralJ1}
	\int_{-1}^1 \frac{J_1\big(\beta\, \sqrt{w^2-1}\big)}{\sqrt{w^2-1}}\,{\mathrm d}w
	&= \int_{-1}^1 \frac{I_1\big(\beta\, \sqrt{1-w^2}\big)}{\sqrt{1-w^2}}\,{\mathrm d}w = \int_{-\pi/2}^{\pi/2} I_1(\beta\,\cos s)\,{\mathrm d}s \notag \\
	&= \pi\,\left(I_{1/2}\Big(\frac{\beta}{2}\Big)\right)^2
	= \frac{4}{\beta}\,\left(\sinh \frac{\beta}{2}\right)^2 \,.
\end{align}
Then from \eqref{eq:eta(v)} and \eqref{eq:integralJ1} it follows that
\begin{align}
\label{eq:eta(v)_case1}
	\eta(v) &= \frac{\beta}{4\, \big(\sinh \frac{\beta}{2}\big)^2}\,\int_{-1}^1 \frac{I_1\big(\beta\, \sqrt{1 -w^2}\big)}{\sqrt{1 -w^2}}\,{\mathrm d}w - \frac{\beta}{2 \, \sinh \beta}\, \int_{-w_1(-v)}^{w_1(v)} \frac{I_1\big(\beta\, \sqrt{1 -w^2}\big)}{\sqrt{1 - w^2}}\,{\mathrm d}w \\[1ex]
	&= \eta_1(v) + \eta_2(v) \notag
\end{align}
with
\begin{eqnarray*}
	\eta_1(v) &\coloneqq& \bigg(\frac{\beta}{4\, \big(\sinh \frac{\beta}{2}\big)^2} - \frac{\beta}{2 \, \sinh \beta}\bigg)\, \int_{-w_1(-v)}^{w_1(v)} \frac{I_1\big(\beta\, \sqrt{1 - w^2}\big)}{\sqrt{1 - w^2}}\,{\mathrm d}w \,, \\ [1ex]
	\eta_2(v) &\coloneqq& \frac{\beta}{4\, \big(\sinh \frac{\beta}{2}\big)^2}\,\bigg(\int_{-1}^1 - \int_{-w_1(-v)}^{w_1(v)}\bigg)\, \frac{I_1\big(\beta\, \sqrt{1 -w^2}\big)}{\sqrt{1 -w^2}}\,{\mathrm d}w \,.
\end{eqnarray*}
By $2\,\big(\sinh \frac{\beta}{2}\big)^2 < \sinh \beta$ we have
\begin{align}
\label{eq:const}
	\frac{\beta}{4\, \big(\sinh \frac{\beta}{2}\big)^2} - \frac{\beta}{2 \, \sinh \beta} > 0\,.
\end{align}
Since the integrand of \eqref{eq:eta(v)_case1} is positive, it is easy to find an upper bound of $\eta_1(v)$ for all $v \in [-\delta,\,\delta]$, because by \eqref{eq:integralJ1} it holds
\begin{eqnarray*}
	0 \le \eta_1(v) &\le& \bigg(\frac{\beta}{4\, \big(\sinh \frac{\beta}{2}\big)^2} - \frac{\beta}{2 \, \sinh \beta}\bigg)\, \int_{-1}^1 \frac{I_1\big(\beta\, \sqrt{1 - w^2}\big)}{\sqrt{1 -w^2}}\,{\mathrm d}w \nonumber\\ [1ex]
	&=& 1 - \frac{2\,\big(\sinh \frac{\beta}{2}\big)^2}{\sinh \beta} = \frac{2 - 2\,{\mathrm e}^{-\beta}}{{\mathrm e}^{\beta} - {\mathrm e}^{-\beta}} < \frac{2}{1 - {\mathrm e}^{-2 \beta}}\, {\mathrm e}^{-\beta}\,.
\end{eqnarray*}
Further, for arbitrary $v \in [- \delta, \, \delta]$ we obtain
\begin{eqnarray}
\label{eq:esteta2}
	0 \le \eta_2(v) &=& \frac{\beta}{4 \, \big(\sinh \frac{\beta}{2}\big)^2}\,\bigg(\int_{-1}^{-w_1(-v)} + \int_{w_1(v)}^1\bigg) \,\frac{I_1\big(\beta\, \sqrt{1 - w^2}\big)}{\sqrt{1 - w^2}}\,{\mathrm d}w \nonumber \\
	&=&  \frac{\beta}{4 \, \big(\sinh \frac{\beta}{2}\big)^2}\,\bigg(\int_{w_1(-v)}^1 + \int_{w_1(v)}^1\bigg) \,\frac{I_1\big(\beta\, \sqrt{1 - w^2}\big)}{\sqrt{1 - w^2}}\,{\mathrm d}w \nonumber \\
	&\le& \frac{\beta}{2 \, \big(\sinh \frac{\beta}{2}\big)^2}\, \int_{w_0}^1 \frac{I_1\big(\beta\, \sqrt{1 - w^2}\big)}{\sqrt{1 - w^2}}\,{\mathrm d}w\,,
\end{eqnarray}
since the integrand is positive and $w_0 \coloneqq w_1(-\delta) = \min \{w_1(v):\,v\in [-\delta,\,\delta]\}$. Substituting $w = \sin t$ in \eqref{eq:esteta2}, for all $v \in [-\delta,\,\delta]$ we can estimate
$$
\eta_2(v) \le \frac{\beta}{2 \, \big(\sinh \frac{\beta}{2}\big)^2}\, \int_{\mathrm{arcsin}\,w_0}^{\pi/2} I_1(\beta\, \cos t)\,{\mathrm d}t
$$
with $\mathrm{arcsin}\,w_0 \in \big(0,\, \frac{\pi}{2}\big)$.
The above integral can now be approximated by the rectangular rule (see Figure~\ref{fig:integrand_i1}) such that
\begin{align*}
	\eta_2(v) \le \frac{\beta}{2\, \big(\sinh \frac{\beta}{2}\big)^2}\,\Big(\frac{\pi}{2} - {\mathrm{arcsin}}\,w_0\Big)\,
	I_1\Big(\beta \,\sqrt{1 - w_0^2}\,\Big) \,.
\end{align*}
Further it holds
$4\, (\sinh \frac{\beta}{2})^2 = {\mathrm e}^{\beta} - 2 + {\mathrm e}^{-\beta} > {\mathrm e}^{\beta} - 2$.
Since by \cite[Lemma 7]{PT21a} we have $\sqrt{2 \pi x}\,{\mathrm e}^{-x}\,I_1(x) < 1$, it holds
$$
I_1\Big(\beta \,\sqrt{1 - w_0^2}\,\Big) < \frac{1}{\sqrt{2 \pi\beta}}\,\big(1 - w_0^2\big)^{-1/4}\,{\mathrm e}^{\beta\,\sqrt{1 - w_0^2}}\,,
$$
and therefore we obtain
\begin{equation*}
	\eta_2(v) \le \frac{\sqrt{\beta}\,(\pi - 2\,{\mathrm{arcsin}}\,w_0)}{\sqrt{2 \pi}\,(1 - w_0^2)^{1/4}\,(1 - 2\,{\mathrm e}^{-\beta})}\,{\mathrm e}^{-\beta\,\big(1 - \sqrt{1-w_0^2}\big)} \,.
\end{equation*}
Additionally using \eqref{eq:E1} and \eqref{eq:eta} as well as $\mathrm{arcsin}\,w_0 \in \big(0,\, \frac{\pi}{2}\big)$ this yields
\begin{align}
\label{eq:err_const_case1}
	E_1(m,\delta,L) \le \frac{\sqrt{\beta\,\pi\delta}}{(1 - 2\,{\mathrm e}^{-\beta})(1-w_0^2)^{1/4}}\,{\mathrm e}^{-\beta\,\big(1-\sqrt{1-w_0^2}\big)} + \frac{2\,\sqrt{2 \delta}}{1 - {\mathrm e}^{-2\beta}}\,{\mathrm e}^{-\beta}\,.
\end{align}

What remains is the choice of the optimal parameter \mbox{$s > 0$}, where we have to fulfill \mbox{$w_1(\delta) \leq 1$.}
To obtain the smallest possible error bound we are looking for an $s>0$ that minimizes the error term $\max_{v\in [-\delta,\,\delta]} |\eta(v)|$.
By \eqref{eq:eta(v)_case1} and \eqref{eq:const} we maximize the second integral in \eqref{eq:eta(v)_case1}.
Since the integrand of \eqref{eq:eta(v)_case1} is positive, the integration limit $w_1(v)$ should be as large as possible for all \mbox{$v \in [-\delta,\, \delta]$} and therefore $w_1(\delta) = 1$.
Rearranging this by \eqref{eq:w1(v)} in terms of~$s$ we see immediately that
$$
s = \frac{m\,(1+\lambda+2\tau)}{1+2\lambda}
$$
and hence
$$
\beta = \frac{\pi m\,(1 + \lambda + 2 \tau)}{1 + \lambda}\,, \quad w_0 = w_1(-\delta) =  \frac{1+\lambda-2\tau}{1+\lambda+2\tau} \in (0,\,1)
$$
such that $\beta$ depends linearly on $m$ by definition.

\begin{figure}[t]
	\centering
	\captionsetup[subfigure]{justification=centering}
	\begin{subfigure}[t]{0.4\textwidth}
		\includegraphics[width=\textwidth]{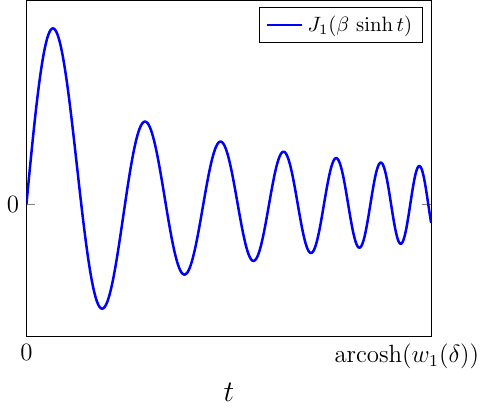}
	\end{subfigure}
	\caption{The integrand $J_1(\beta\sinh t)$ on the interval $[0,\mathrm{arcosh}(w_1(\delta))]$.}
	\label{fig:integrand_j1}
\end{figure}

\paragraph{Case 2 ($\boldsymbol{w_1(-\delta)\geq 1}$):}
From \eqref{eq:eta(v)} it follows that
\begin{align}
\label{eq:eta(v)_case2}
	\eta(v) = \eta_3(v) - \eta_4(v)\,, \quad v \in [-\delta,\, \delta]\,,
\end{align}
with
\begin{eqnarray*}
	\eta_3(v) &\coloneqq& 1 - \frac{\beta}{2\,\sinh \beta}\,\int_{-1}^1 \frac{I_1\big(\beta\, \sqrt{1-w^2}\big)}{\sqrt{1-w^2}}\, {\mathrm d}w\,, \\ [1ex]
	\eta_4(v) &\coloneqq& \frac{\beta}{2\,\sinh \beta}\,\bigg(\int_{-w_1(-v)}^{-1} + \int_1^{w_1(v)}\bigg) \,\frac{J_1\big(\beta\, \sqrt{w^2-1}\big)}{\sqrt{w^2-1}}\, {\mathrm d}w\,.
\end{eqnarray*}
By \eqref{eq:integralJ1} we obtain
$$
\eta_3(v) = 1 - \frac{2 \big(\sinh \frac{\beta}{2}\big)^2}{\sinh\,\beta} = \frac{2\, {\mathrm e}^{-\beta}}{1 + {\mathrm e}^{-\beta}} > 0\,.
$$
Further it holds
$$
\eta_4(v) = \frac{\beta}{2\,\sinh \beta}\,\bigg(\int_1^{w_1(-v)} + \int_1^{w_1(v)}\bigg) \frac{J_1\big(\beta\, \sqrt{w^2-1}\big)}{\sqrt{w^2-1}}\, {\mathrm d}w\,.
$$
Substituting $w=\cosh t$ in above integrals, we have
$$
\eta_4(v) = \frac{\beta}{2\,\sinh \beta}\,\bigg( \int_0^{\mathrm{arcosh}(w_{1}(-v))} +\int_0^{\mathrm{arcosh}(w_{1}(v))} \bigg) J_1(\beta\, \sinh t) \,{\mathrm d}t \,.
$$
In order to estimate these integrals properly, we now have a closer look at the integrand. As known,
the Bessel function $J_1$ oscillates on $[0,\,\infty)$ and has the non-negative simple zeros $j_{1,n}$, $n \in {\mathbb N}_0$, with $j_{1,0} = 0$. The zeros $j_{1,n}$, $n = 1,\,\ldots,\,40$, are tabulated in \cite[p.~748]{Wa44}.
On each interval $\big[\mathrm{arsinh}\,\frac{j_{1,2n}}{\beta},\,\mathrm{arsinh}\,\frac{j_{1,2n+2}}{\beta}\big]$, $n \in {\mathbb N}_0$, the integrand $J_1(\beta\,\sinh t)$ is firstly non-negative and then non-positive, see Figure~\ref{fig:integrand_j1}.
Due to this properties and the fact that the amplitude is decreasing when $x\to\infty$, the integrals are positive on each interval $\big[\mathrm{arsinh}\,\frac{j_{1,2n}}{\beta},\,\mathrm{arsinh}\,\frac{j_{1,2n+2}}{\beta}\big]$, $n \in {\mathbb N}_0$.
Note that by \cite[6.645--1]{GR80} it holds
$$\int_0^{\infty} J_1(\beta \sinh t) \,{\mathrm d}t = I_{1/2}\Big(\frac{\beta}{2}\Big) \, K_{1/2}\Big(\frac{\beta}{2}\Big) = \frac{2}{\sqrt{\pi\beta}}\,\sinh \frac{\beta}{2} \cdot \sqrt{\frac{\pi}{\beta}}\,\e^{-\beta/2} = \frac{1 - \e^{- \beta}}{\beta} \,,$$
where $K_{\alpha}$ denotes the modified Bessel function of second kind and $I_{1/2}$, $K_{1/2}$ denote modified Bessel functions of half order (see \cite[10.2.13, 10.2.14, and 10.2.17]{abst}).
In addition, numerical experiments have shown that for all $T\ge 0$ it holds
$$
0\leq\int_0^{T} J_1(\beta \sinh t) \,{\mathrm d}t\leq \frac{3\,(1 - \e^{- \beta})}{2 \, \beta} \,.
$$
Therefore, we obtain
$$
0 \le \eta_4(v) \le \frac{\beta}{2\, \sinh \beta} \cdot \frac{3\,(1 - {\mathrm e}^{-\beta})}{\beta} = \frac{3\,{\mathrm e}^{-\beta}}{1 + {\mathrm e}^{-\beta}}< 3\,{\mathrm e}^{-\beta}
$$
and hence by \eqref{eq:eta(v)_case2} it holds
\begin{align}
\label{eq:est_eta_case2}
	\max_{v\in [-\delta,\,\delta]} |\eta(v)| = \max_{v\in [-\delta,\,\delta]} |\eta_3(v) - \eta_4(v)| < 3\,{\mathrm e}^{-\beta}\,.
\end{align}
Thus, by \eqref{eq:E1} and \eqref{eq:eta} we conclude that
\begin{align}
\label{eq:err_const_case2}
	E_1(m,\delta,L)
	\le
	3\,\sqrt{2\delta}\, \e^{-\beta} \,.
\end{align}

What remains is the choice of the optimal parameter \mbox{$s > 0$}, where we have to fulfill \mbox{$w_1(-\delta) = c$} with $c\geq 1$.
Rearranging this by \eqref{eq:w1(v)} in terms of $s$ we see that
$$
s = s(c) = \frac{m\,(1 + \lambda - 2 \tau)}{c\,(1 + 2\lambda)}\,, \quad \beta = \beta(c) = \frac{\pi m\,(1 + \lambda - 2 \tau)}{c\,(1 + \lambda)}.
$$
To obtain the smallest possible error bound we are looking for a constant $c\geq 1$ that minimizes the error term $\max_{v\in [-\delta,\,\delta]} |\eta(v)|$.
By \eqref{eq:est_eta_case2} we minimize the upper bound
$3\,{\mathrm e}^{-\beta(c)}$.
Since $3\,{\mathrm e}^{-\beta(c)}$ is monotonously increasing for $c\geq 1$, we recognize that the minimum value is obtained for $c=1$.
Hence, the suggested parameters are
$$
s = \frac{m\,(1 + \lambda - 2 \tau)}{1 + 2\lambda}\,, \quad \beta = \frac{\pi m\,(1 + \lambda - 2 \tau)}{1 + \lambda}
$$
such that $\beta$ depends linearly on $m$ by definition.
This completes the proof.
\qedsymbol
\medskip

Now we compare the actual decay rates of the error constants \eqref{eq:err_const_case1} with $\beta = \frac{\pi m\,(1 + \lambda +2 \tau)}{1 + \lambda}$ and \eqref{eq:err_const_case2} with $\beta = \frac{\pi m\,(1 + \lambda - 2 \tau)}{1 + \lambda}$.
It can be seen that the decay rate of \eqref{eq:err_const_case1} reads as
$$
\frac{\pi m\,(1 + \lambda +2 \tau)}{1 + \lambda}\,\Big(1 - \sqrt{1 - w_0^2}\,\Big)
$$
with $w_0 = \frac{1 + \lambda - 2\tau}{1 + \lambda + 2\tau}$.
On the other hand, the decay rate of \eqref{eq:err_const_case2} is given by $\frac{\pi m\,(1 + \lambda - 2 \tau)}{1 + \lambda}$. Since \mbox{$1 + \lambda > 2\tau$} for all $\lambda \ge 0$ and $\tau \in \big(0,\,\frac{1}{2}\big)$, simple calculation shows that
$$
\frac{\pi m\,(1 + \lambda - 2 \tau)}{1 + \lambda} > \frac{\pi m\,(1 + \lambda +2 \tau)}{1 + \lambda}\,\Big(1 - \sqrt{1 - w_0^2}\,\Big)\,.
$$
Hence, the error constant \eqref{eq:err_const_case2} decays faster than the one in \eqref{eq:err_const_case1}.
Therefore, we will use the $\sinh$-type window function \eqref{eq:varphisinh} with $\beta = \frac{\pi m\,(1 + \lambda - 2 \tau)}{1 + \lambda}$ in the remainder of this paper.

\begin{Example}
\label{ex:visualization_sinhWKS}
Analogous to Example~\ref{ex:visualization_GaussianWKS}, we visualize the error bound of Theorem \ref{Thm:errorpsisinh}, i.\,e.,
for $\varphi = \varphi_{\mathrm{sinh}}$ we show that for the approximation error \eqref{eq:e_ndelta} it holds by \eqref{eq:err_const_case2} that $e_{m,\tau,\lambda}(f) \le E_{m,\tau,\lambda}\,\|f\|_{L^2(\mathbb R)}$, where
\begin{equation}
	\label{eq:E_ndelta_sinh}
	E_1(m,\delta,L)
	\le  E_{m,\tau,\lambda} \coloneqq
	3\,\sqrt{2\delta} \, \e^{-\beta} \,,
\end{equation}
with \mbox{$\beta = \frac{\pi m\,(1+\lambda - 2 \tau)}{1+\lambda}$}.
The associated results are displayed in Figure~\ref{fig:err_const_sinh}, where a substantial improvement can be seen compared to the Figures~\ref{fig:err_const_Gaussian} and \ref{fig:err_const_bspline}.
We also remark that for larger choices of $N$, the line plots in Figure~\ref{fig:err_const_sinh} would only be shifted slightly upwards, such that for all $N$ we receive almost the same results.
This is to say, the $\sinh$-type window function is by far the best choice as a regularization function for regularized Shannon sampling sums.
\begin{figure}[ht]
	\centering
	\captionsetup[subfigure]{justification=centering}
	\begin{subfigure}[t]{0.49\textwidth}
		\includegraphics[width=\textwidth]{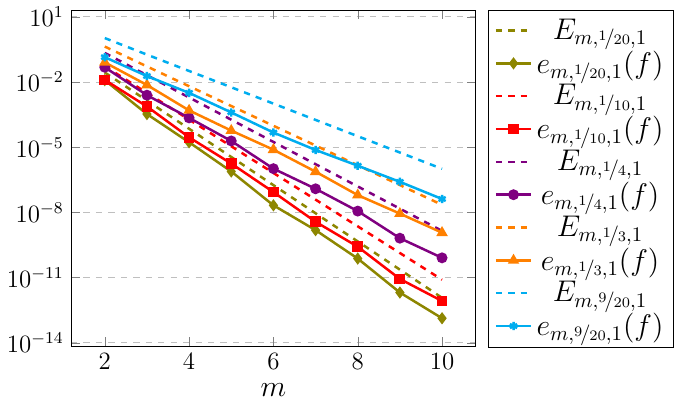}
		\caption{{$\lambda=1$} and various $\tau<\frac 12$}
	\end{subfigure}
	\begin{subfigure}[t]{0.49\textwidth}
		\includegraphics[width=\textwidth]{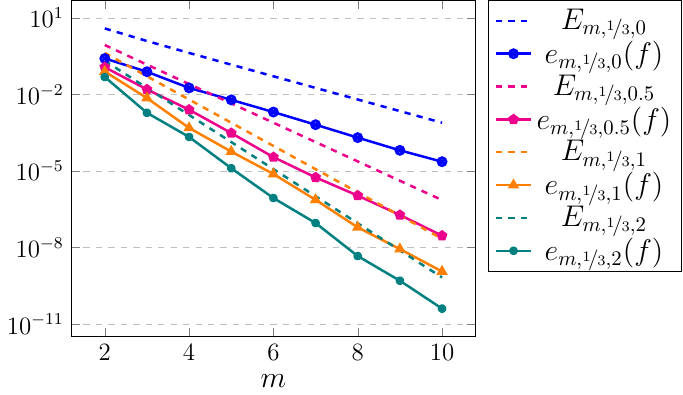}
		\caption{$\tau=\frac 13$ and various $\lambda\geq 0$}
	\end{subfigure}
	\caption{Maximum approximation error \eqref{eq:e_ndelta} and error constant \eqref{eq:E_ndelta_sinh} using \mbox{$\varphi_{\sinh}$} and \mbox{$\beta = \frac{\pi m\,(1 + \lambda - 2 \tau)}{1 + \lambda}$} in \eqref{eq:varphisinh} for the function \mbox{$f(x) = \sqrt{2\delta} \,\mathrm{sinc}(2 \delta \pi x)$} with \mbox{$N=128$}, \mbox{$m\in\{2, 3, \ldots, 10\}$}, as well as \mbox{$\tau \in \{\nicefrac{1}{20},\,\nicefrac{1}{10},\,\nicefrac{1}{4},\,\nicefrac{1}{3},\,\nicefrac{9}{20}\}$}, \mbox{$\delta = \tau N$}, and \mbox{$\lambda\in\{0,0.5,1,2\}$}, respectively.
		\label{fig:err_const_sinh}}
\end{figure}
\end{Example}

Now we show that for the regularized Shannon sampling formula with the $\sinh$-type window function \eqref{eq:varphisinh} the uniform perturbation error \eqref{eq:result_robustness}
only grows as $\O{\sqrt{m}}$.

\begin{Theorem}
\label{Theorem:robustness_sinh}
Let $f\in {\mathcal B}_{\delta}(\mathbb R)$ with $\delta = \tau N$, $\tau\in (0,\,1/2)$, $N \in \mathbb N$, \mbox{$L= N(1+\lambda)$} with \mbox{$\lambda\geq 0$} and \mbox{$m \in {\mathbb N}\setminus \{1\}$} be given.
Further let \mbox{$R_{\sinh,m}{\tilde f}$} be as in \eqref{eq:Rmf_erroneous} with the noisy samples
${\tilde f}_{\ell} = f\big(\frac{\ell}{L}\big) + \varepsilon_{\ell}$,
where $|\varepsilon_{\ell}| \le \varepsilon$ for all $\ell\in\Z$ and $\varepsilon>0$.\\
Then the regularized Shannon sampling formula \eqref{eq:Rmf} with the $\sinh$-type window function~\eqref{eq:varphisinh} and $\beta = \frac{\pi m\,(1 + \lambda - 2 \tau)}{1 + \lambda}$ satisfies
\begin{equation}
\label{eq:result_robustness_sinh}
	\| R_{\sinh,m}{\tilde f} - R_{\sinh,m}f \|_{C_0(\mathbb R)}
	\leq \varepsilon \left( 2+\sqrt{\frac{2+2\lambda}{1+\lambda - 2 \tau}} \,\frac{1}{1 - {\mathrm e}^{-2 \beta}}\,\sqrt{m} \right)\,.
\end{equation}
\end{Theorem}

{\emph Proof.} By Theorem \ref{Theorem:robustness} we only have to compute $\hat{\varphi}_{\sinh}(0)$.
By \eqref{eq:sinh_FT} and $\sqrt{2 \pi \beta}\,{\mathrm e}^{-\beta}\,I_1(\beta) < 1$ (see \cite[Lemma 7]{PT21a}) we recognize that
\begin{align*}
\hat{\varphi}_{\mathrm{sinh}}(0)
&=
\frac{\pi m \beta}{L\, \sinh \beta} \cdot \frac{I_1\big(\beta\big)}{\beta}
=
\frac{\pi m\,I_1\big(\beta\big)}{L\, \sinh \beta}
\le
\frac{\pi m \,{\mathrm e}^{\beta}}{\sqrt{2 \pi \beta}L\, \sinh \beta}=\frac{\sqrt{2\pi}\,m}{\sqrt{\beta}L\,(1 - {\mathrm e}^{- 2 \beta})}\,.
\end{align*}
If we now use $\beta = \frac{\pi m\,(1+\lambda - 2 \tau)}{1+\lambda}$, then
\eqref{eq:result_robustness} yields the assertion \eqref{eq:result_robustness_sinh}.
\qedsymbol
\medskip

Similar to Example~\ref{ex:perturbation_Gaussian}, one can visualize the error bound from Theorem~\ref{Theorem:robustness_sinh}, which leads to results analogous to Figure~\ref{fig:err_perturbation_Gaussian}.

\section{Conclusion\label{sec:Conclusion}}

To overcome the drawbacks of classical Shannon sampling series -- which are poor convergence and non-robustness in the presence of noise -- in this paper we considered regularized Shannon sampling formulas with localized sampling.
To this end, we considered bandlimited functions $f\in\mathcal{B}_{\delta}(\R)$ and introduced a set $\Phi_{m,L}$ of window functions.
Despite the original result, where $\varphi \in \Phi_{m,L}$ is chosen as the rectangular window function, and the well--studied approach of using the Gaussian window function, we proposed new window functions with compact support $[-m/L,\,m/L]$, namely the $\mathrm B$--spline and $\sinh$-type window function, which are well-studied in the context of the nonequispaced fast Fourier transform (NFFT).

In Section~\ref{sec:WKS_localized} we considered an arbitrary window function $\varphi \in \Phi_{m,L}$ and presented a unified approach to error estimates of the uniform approximation error in Theorem~\ref{Theorem:errorRmf}, as well as a unified approach to the numerical robustness in Theorem~\ref{Theorem:robustness}.

Then, in the next sections, we concretized the results for special choices of the window functions.
More precisely, it was shown that the uniform approximation error decays exponentially with respect to the truncation parameter $m$, if $\varphi \in \Phi_{m,L}$ is the Gaussian, $\mathrm B$--spline, or $\sinh$-type window function.
Moreover, we have shown that the regularized Shannon sampling formulas are numerically robust for noisy samples, i.\,e., if $\varphi \in \Phi_{m,L}$ is the Gaussian, $\mathrm B$--spline, or $\sinh$-type window function, then the uniform perturbation error only grows as $m^{1/2}$.
While the Gaussian window function from Section~\ref{sec:Gaussian} has already been studied in numerous papers such as \cite{Q03, Q04, Q05, Q06, LZ16}, we remarked that Theorem~\ref{Theorem:GaussianWKS} improves a corresponding result in \cite{LZ16}, since we improved the exponential decay rate from $(m-1)$ to $m$.

\begin{figure}[ht]
	\centering
	\captionsetup[subfigure]{justification=centering}
	\begin{subfigure}[t]{0.32\textwidth}
		\includegraphics[width=\textwidth]{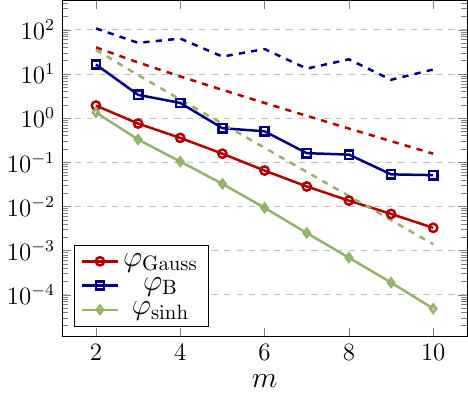}
		\caption{$\lambda=0.5$}
	\end{subfigure}
	\begin{subfigure}[t]{0.32\textwidth}
		\includegraphics[width=\textwidth]{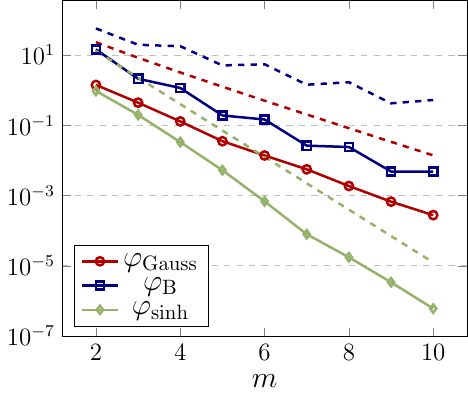}
		\caption{$\lambda=1$}
	\end{subfigure}
	\begin{subfigure}[t]{0.32\textwidth}
		\includegraphics[width=\textwidth]{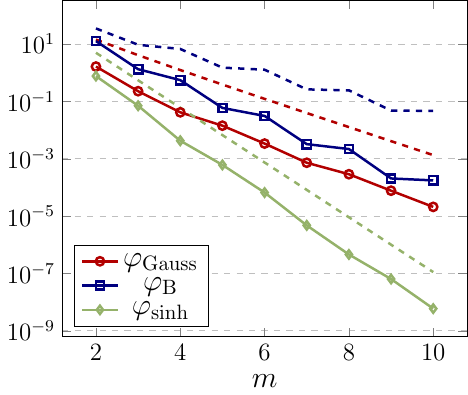}
		\caption{$\lambda=2$}
	\end{subfigure}
	\caption{Maximum approximation error \eqref{eq:e_ndelta} and error constant \eqref{eq:error_const} using \mbox{$\varphi\in\{\varphi_{\mathrm{Gauss}},\varphi_{\mathrm{B}}, \varphi_{\sinh}\}$} for the function \mbox{$f(x) = \delta \,\mathrm{sinc}^2(\delta \pi x)$} with \mbox{$N=256$}, \mbox{$\tau = 0.45$}, \mbox{$\delta = \tau N$}, as well as \mbox{$m\in\{2, 3, \ldots, 10\},$} and \mbox{$\lambda\in\{0.5,1,2\}$}.}
	\label{fig:comp_err}
\end{figure}

Throughout this paper, several numerical experiments illustrated the corresponding theo\-retical results.
Finally, comparing the proposed window functions as done in Figure~\ref{fig:comp_err}, the superiority of the new proposed $\sinh$-type window function can easily be seen, since even small choices of the truncation parameter $m\leq 10$ are sufficient for achieving high precision.
Due to the usage of localized sampling the evaluation of $R_{\varphi,m}f$ on an interval $[0,1/L]$ requires only $2m$ samples and therefore has a computational cost of $\O{2m}$ flops.
Thus, a reduction of the truncation parameter $m$ is desirable to obtain an efficient method.

\section*{Acknowledgments}
Melanie Kircheis gratefully acknowledges the support from the BMBF grant 01$\mid$S20053A (project SA$\ell$E).
Daniel Potts acknowledges the funding by Deutsche Forschungsgemeinschaft (German Research Foundation) -- Project--ID 416228727 -- SFB 1410.

Moreover, the authors thank the referees and the editor for their very useful suggestions for improvements.


\begin{thebibliography}{12}

\bibitem{abst}
M.~Abramowitz and I.A.~Stegun, editors.
\newblock {\em Handbook of Mathematical Functions with Formulas, Graphs, and Mathematical Tables}.
\newblock Dover, New York, 1972.

\bibitem{chui92}
C.~K.~Chui.
\newblock {\em An Introduction to Wavelets}.
\newblock {Academic Press, Boston, 1992}.

\bibitem{Const16}
A.~Constantin.
\newblock {\em Fourier Analysis, Part I: Theory}.
\newblock Cambridge University Press, Cambridge, 2016.

\bibitem{DDeV03}
I.~Daubechies and R.~DeVore.
\newblock{Approximating a bandlimited function using very coarsely quantized data: A family of stable sigma-delta modulators of arbitrary order}.
\newblock{\em Ann. of Math. (2)}, 158:679--710, 2003.

\bibitem{GR80}
I.S.~Gradshteyn and I.M.~Ryzhik.
\newblock{\em Table of Integrals, Series, and Products}.
\newblock Academic Press, New York, 1980.

\bibitem{LZ16}
R.~Lin and H.~Zhang.
\newblock {Convergence analysis of the Gaussian regularized Shannon sampling formula}.
\newblock{\em Numer. Funct. Anal. Optim.}, 38(2):224--247, 2017.

\bibitem{MeRo65}
R.G.~Medhurst and J.H.~Roberts.
\newblock{Evaluation of the integral $I_n(b) = \frac{2}{\pi}\,\int_0^{\infty} (\sin\,x/x)^n\,\cos(b x)\,{\mathrm d}x$}.
\newblock{\em Math. Comp.}, 19:113--117, 1965.

\bibitem{MXZ09}
C.A.~Micchelli, Y.~Xu, and H.~Zhang.
\newblock {Optimal learning of bandlimited functions from localized sampling}.
\newblock {\em J. Complexity}, 25(2):85--114, 2009.

\bibitem{Ob90}
F.~Oberhettinger.
\newblock {\em Tables of Fourier Transforms and Fourier Transforms of Distributions}.
\newblock Springer, Berlin, 1990.

\bibitem{PPST18}
G.~Plonka, D.~Potts, G.~Steidl, and M.~Tasche.
\newblock {\em Numerical Fourier Analysis}.
\newblock Birkh\"auser/Springer, Cham, 2018.

\bibitem{PT21a}
D.~Potts and M.~Tasche.
\newblock{Uniform error estimates for nonequispaced fast Fourier transforms}.
\newblock{\em  Sampl. Theory Signal Process. Data Anal.}, 2021.

\bibitem{PT21}
D.~Potts and M.~Tasche.
\newblock {Continuous window functions for NFFT}.
\newblock {\em Adv. Comput. Math.}, 2021.

\bibitem{Q03}
L.~Qian.
\newblock {On the regularized Whittaker--Kotelnikov--Shannon sampling formula}.
\newblock {\em Proc. Amer. Math. Soc.}, 131(4):1169--1176, 2003.

\bibitem{Q04}
L.~Qian.
\newblock {The regularized Whittaker-Kotelnikov-Shannon sampling theorem and its application to the numerical solutions of partial differential equations}.
\newblock{\em PhD thesis, National Univ. Singapore}, 2004.

\bibitem{Q05}
L.~Qian and D.B.~Creamer.
\newblock {Localized sampling in the presence of noise}.
\newblock {\em Appl. Math. Letter}, 19:351--355, 2006.

\bibitem{Q06}
L.~Qian and D.B.~Creamer.
\newblock {A modification of the sampling series with a Gaussian multiplier}.
\newblock {\em Sampl. Theory Signal Image Process.}, 5(1):1--20, 2006.

\bibitem{SchSt07}
G.~Schmeisser and F.~Stenger.
\newblock {Sinc approximation with a Gaussian multiplier}.
\newblock {\em Sampl. Theory Signal	Image Process.}, 6(2):199--221, 2007.

\bibitem{StTa06}
T.~Strohmer and J.~Tanner.
\newblock {Fast reconstruction methods for bandlimited functions from periodic nonuniform sampling}.
\newblock {\em SIAM J. Numer. Anal.}, 44(3):1071-1094, 2006.

\bibitem{TaSuMu07}
K.~Tanaka, M.~Sugihara, and K.~Murota.
\newblock {Complex analytic approach to the sinc-Gauss sampling formula}.
\newblock {\em Japan J. Ind. Appl. Math.}, 25:209--231, 2008.

\bibitem{UAE}
M.~Unser, A.~Aldroubi, and M.~Eden.
\newblock {On the asymptotic convergence of B-spline wavelets to Gabor functions}.
\newblock {\em IEEE Trans. Inform. Theory} 38(2):864--872, 1992.

\bibitem{Wa44}
G.N.~Watson.
\newblock {\em A Treatise on the Theory of Bessel Functions}, 2nd edn.
\newblock Cambridge Univ. Press, Cambridge, 1944.

\end{thebibliography}
\end{document}